\documentclass[12pt]{amsart}

\usepackage{epsfig}
\usepackage{subfig}
\usepackage{amscd}
\usepackage[mathscr]{eucal}
\usepackage{amssymb}
\usepackage{amsxtra}
\usepackage{amsmath}
\usepackage[all]{xy}
\usepackage{tikz}
\usepackage{pgf}
\usetikzlibrary{calc,backgrounds,arrows,matrix,shapes}
\usepackage{verbatim}

\theoremstyle{plain}

\theoremstyle{definition}

\theoremstyle{remark}

\begin{document}

\bigskip

\title[E-collections and mirror symmetry for toric Fano]{On exceptional collections of line bundles and mirror symmetry for toric Del-Pezzo surfaces}

\email{yochay.jerby@gmail.com}

\date{\today}

\author{Yochay Jerby}


%
%
\begin{abstract}
Let $X$ be a toric Del-Pezzo surface and let $Crit(W) \subset (\mathbb{C}^{\ast})^n$ be the solution scheme of the Landau-Ginzburg system of equations. Denote by $X^{\circ}$ the polar variety of $X$. Our aim in this work is to describe a map $L : Crit(W) \rightarrow Fuk_{trop}(X^{\circ})$
whose image under homological mirror symmetry corresponds to a full strongly exceptional collection of line bundles.

\end{abstract}

\maketitle

%
%
\section{Introduction and Summary of Main Results}
\label{s:intro}

\bigskip
\hspace{-0.6cm} Let $X$ be a smooth toric Fano manifold given by a reflexive polytope $ \Delta \subset \mathbb{R}^n$ and let $\Delta^{\circ} \subset (\mathbb{R}^n)^{\vee}$ be the polar polytope of $\Delta$, see \cite{Ba2}. The polytope $\Delta^{\circ}$ determines $\widetilde{L}(\Delta^{\circ})$, the family of Laurent polynomials $ W : (\mathbb{C}^{\ast})^n \rightarrow \mathbb{C}$ whose Newton polytope is $\Delta^{\circ}$. We refer to a pair $((\mathbb{C}^{\ast})^n,W)$, where $W \in \widetilde{L}(\Delta^{\circ})$ is such that $\left \{ W=0 \right \}$ is smooth, as a mirror model of $X$.

\hspace{-0.6cm} In \cite{Ab1,Ab2} Abouzaid established an equivalence between $\mathcal{D}^b(X)$, the derived category of coherent sheaves of $X$, and a certain full subcategory of $\mathcal{D}^{\pi} (Fuk( (\mathbb{C}^{\ast})^n, \widetilde{W}))$, the derived Fukaya category of a specific mirror model $((\mathbb{C}^{\ast})^n, \widetilde{W})$, affirming Kontsevich's homological mirror symmetry conjecture for toric Fano manifolds, see \cite{Ko,Ko2}. The equivalence is based on the definition of
a correspondence between the $T$-equivariant line bundles $\mathcal{O}(D)$, for $D \in Div_T(X)$, and hamiltonian isotopy classes of certain Lagrangian branes in $Fuk((\mathbb{C}^{\ast})^n, \widetilde{W})$, called tropical sections. The extension to the whole of $\mathcal{D}^b(X)$ relies on the homolgical algebraic fact that, in the toric case, $\mathcal{D}^b(X)$ is generated as a triangulated category by the line bundles $\mathcal{O}(D)$, for all $ D \in Div_T(X)$. We hence-forward denote by $Fuk(X^{\circ})$ for the Fukaya pre-category of the model $((\mathbb{C}^{\ast})^n,\widetilde{W})$ and by $Fuk_{trop}(X^{\circ})$ the sub-pre-category of tropical sections.

\hspace{-0.6cm} One of the features of Abouzaid's homological mirror symmetry functor is that it does not make reference to exceptional collections on neither of the categories $ \mathcal{D}^b(X)$ and $\mathcal{D}^{\pi}( Fuk(X^{\circ}))$, for its definition. On the other hand, one of the general approaches to establish homological mirror symmetry is to identify collections $\mathcal{L}$ of Lagrangian submanifolds in the mirror which are analogs of known full strongly exceptional collections $\mathcal{E}$ in $\mathcal{D}^b(X)$. For instance, Seidel's fundamental construction describes such collections of Lagrangian vanishing cycles in the directed Fuaya-Seidel category, which depend on a certain choice of thimbles, see \cite{S,S2,S3}. In a few cases (toric and non-toric), a correspondence between such collections of Lagrangian vanishing cycles and specific known exceptional collections of sheaves $\mathcal{E}$ in $\mathcal{D}^b(X)$ has been established, see for instance \cite{AKO,AKO2,S2,S4,Ue1}. However, in general, it becomes increasingly hard to determine, in advance, what is the appropriate analog in $\mathcal{D}^b(X)$ of a given vanishing cycle $L \in \mathcal{L}$. We refer the reader to \cite{DKK, GW} where the case of general nef toric stacks is studied. In particular, the established examples show that classes analog to vanishing cycles under mirror symmetry typically do not turn to be classes of line bundles, or even vector bundles, in $\mathcal{D}^b(X)$.

\hspace{-0.6cm} Let $W \in \widetilde{L}(\Delta^{\circ})$ be an element and let $Crit(W) \subset (\mathbb{C}^{\ast})^n$ be the solution scheme of the Landau-Ginzburg system of equations given by
$$ z_i \frac{ \partial}{\partial z_i} W(z_1,...,z_n) = 0 \hspace{0.5cm} \textrm{ for } \hspace{0.25cm} i=1,...,n.$$ In \cite{J,J2} we described, for various toric Fano manifolds $X$ and potentials $W$, maps of the form $ E : Crit(W) \rightarrow Pic(X)$ whose image is a full strongly exceptional collection, we refered to these maps as exceptional maps.
Our aim in this work is to describe maps $L : Crit(W) \rightarrow Fuk_{trop}(X^{\circ})$ which make the following diagram commutative

$$\begin{tikzpicture}[every node/.style={midway}]
\matrix[column sep={6em,between origins}, row sep={3em}] at (0,0) { & \node(C) {$Crit(W)$}; & \\
\node(P) {$Pic(X)$}; & & \node (F) {$\frac{Fuk_{trop}(X^{\circ})}{\sim}$};\\
};
\draw[->] (C) -- (P) node[anchor=east] {$E$};
\draw[->] (C) -- (F) node[anchor=west] {$L$};
\draw[->] (P) -- (F) node[anchor=north] {$HMS$};
\end{tikzpicture}$$ where the map $HMS$, in the bottom row, is Abouzaid's homological mirror symmetry functor restricted to $Pic(X)$. In particular, let us note that as in Seidel's construction, our apporoch is to associate Lagrangian branes to elements of $Crit(W)$. However, in our case the corresponding branes are determined by the combinatorial data of the polytope $\Delta$, and do not involve additional choices.
‏

\hspace{-0.6cm} In order to construct $L$ we first inflate the elements of $Crit(W) \subset (\mathbb{C}^{\ast})^n$ into a collection of embeded Lagragian $n$-balls. This is achieved by considering
the solutions of $$ z_i \frac{ \partial}{\partial z_i} W(z_1,...,z_n) = a_i \hspace{0.5cm} \textrm{ for } \hspace{0.25cm} i=1,...,n,$$ with $a \in \mathbb{R}^n$ such that $\vert a \vert \leq r$ for some $0<r$. We denote by $\widetilde{L}_r(z)$ the resulting
Lagrangian ball whose center is $z \in Crit(W)$. In the examples considered one always has a unique solution $z_0$ such that $E(z_0) = \mathcal{O}_X$. It turns that of all the resulting balls, only $\widetilde{L}_r(z_0)$ is a genuine element of $Fuk(X^{\circ})$. In fact, $\widetilde{L}_r(z_0)$ turns to coincide with the trivial tropical section $L_0 \in Fuk_{trop}(X^{\circ})$, which corresponds to $\mathcal{O}_X$ under $HMS$. We hence want to fix the rest of the balls to be tropical sections.
In order to do so consider the monodromy map defined by standard analytic continuation along a loop $$ M : \pi_1 ( L(\Delta^{\circ}) \setminus R_X ,
W) \rightarrow Aut(Crit(W)), $$ where $R_X$ is the hypersurface of all $W$ for which $Crit(W)$ is non-reduced and $Aut(Crit(W))$ is the permutation group of $Crit(W)$ as a finite set, see \cite{J, J2}. The map $M$ inflates to an analog map for the inflated Lagrangian balls by setting $M(\widetilde{L}_r(z)):=\widetilde{L}_r(M(z))$. Let $z \in Crit(W)$ be any solution and set $$G(z,z_0) := \left \{ [\gamma] \vert M(\gamma)(z)=z_0 \right \} \subset \pi_1 (L(\Delta^{\circ}) \setminus R , W).$$
Our main observation is that a maximal cone $\sigma \in \Sigma(n)$ in the fan $\Sigma$ of $X$, see \cite{F,O}, naturally determines an element $\gamma(z ; \sigma) \in G(z,z_0)$. Note
that an element $\gamma \in G(z,z_0)$ can be associated with a weight $m(\gamma) \in \mathbb{Z}^n$ by lifting the image of the monodromy $M(\gamma)$ under the Argument map to the universal cover
$\mathbb{R}^n$ of $\mathbb{T}^n$. We set $m_W(z ; \sigma):= m(\gamma(z ; \sigma)) \in \mathbb{Z}^n$.

\hspace{-0.6cm} On the other hand, recall that $Div_T(X)$ is isomorphic to the group of continuous piecewise linear functions, supported on $\Sigma$. In particular, a divisor
$D \in Div_T(X)$ and a maximal cone $\sigma \in \Sigma(n)$ are associated with a weight $m_X(D ; \sigma ) \in \mathbb{Z}^n$, viewed as the linear functional obtained by restricting the support function of $D$ to $\sigma$. We show:

\bigskip

\hspace{-0.6cm} \bf Theorem A: \rm Let $X$ be a toric Del Pezzo surface and let $z \in Crit(W)$ be the solution of the corresponding potential $W$. Then $$ m_W(z ; \sigma) =-m_X( D_z ; \sigma ) $$ for any $\sigma \in \Sigma(n)$ where $D_z \in Div_T(X)$ is a $T$-divisor
such that $E(z)=[D_z] \in Pic(X)$.

\bigskip

\hspace{-0.6cm} Let us note that we view the above formula as a manifestation of the mirror symmetry phenomena, as it express a relation between invariants of $\mathcal{D}^b(X)$ and geometric properties of the mirror. Finally, the
collection of monodromies $M( \gamma(z ; \sigma))$ allows us to construct the required sections. Indeed, even though the Lagrangian balls $\widetilde{L}_r(z)$ are not tropical sections, for $ z \neq z_0$, they still satisfiy $Log_r (\widetilde{L}_r(z)) \rightarrow \Delta$, due to standard considerations of tropical geometry. We thus view them as embedded copies of (a smooth approximation of) $\Delta$. To each vertex $x \in \Delta(0)$ we apply the monodromy $M(\gamma(z ; \sigma_z))$ where $\sigma_z$ is the corresponding cone
under the identification $\Sigma(n) \simeq \Delta(0)$. To the rest of the one-dimensional boundary points we apply the monodromies obtained by interpolating the paths $\gamma(z ; \sigma)$ associated to the vertices on the endpoints, in the obvious way. In particular, verifying with the
definition of the map $HMS$ together with Theorem A shows that the resulting tropical sections satisfy the required $HMS(E(z)) = [L_r(z)]$.

\bigskip

\hspace{-0.6cm} The rest of the work is organized as follows: In section 2 we review relevant results on toric Fano manifolds and exceptional collections in $Pic(X)$ as well as the construction of the
maps $E$. In section 3 we review relevant features of Abouzaid's homological mirror symmetry functor. In section 4 we define $L$ for the toric Del-Pezzo surfaces and discuss generalization to
higher dimensional cases and concluding remarks.

\section{Review on exceptional collections in $Pic(X)$ of toric Fano manifolds}
\label{s:intro}

\hspace{-0.6cm} A toric variety is an algebraic variety $X$ containing an algebraic torus $T \simeq (\mathbb{C}^{\ast})^n$ as a dense subset such that the action of $T$ on itself extends to the whole variety, we refer the reader to \cite{F,O}. A compact toric variety $X$ is said to be Fano if its anticanonical class $-K_X$ is Cartier and ample.

\hspace{-0.6cm} Let $N \simeq \mathbb{Z}^n$ be a lattice and let $M = N^{\vee}=Hom(N, \mathbb{Z})$ be the dual lattice. Denote by $N_{\mathbb{R}} = N \otimes \mathbb{R}$ and $M_{\mathbb{R}}=M \otimes \mathbb{R}$ the corresponding vector spaces. Let $ \Delta \subset M_{\mathbb{R}}$ be an integral polytope and let $$ L(\Delta) := \bigoplus_{m \in \Delta \cap M} \mathbb{C} m. $$ The polytope $\Delta $ determines the embedding $i_{\Delta} : (\mathbb{C}^{\ast})^n \rightarrow \mathbb{P}(L(\Delta)^{\vee})$ given by $ z \mapsto [z^m \mid m \in \Delta \cap M] $. The \emph{toric variety} $X_{\Delta} \subset \mathbb{P}(L(\Delta)^{\vee})$ corresponding to the polytope $\Delta \subset M_{\mathbb{R}}$ is defined
to be the compactification of the embedded torus $i_{\Delta}((\mathbb{C}^{\ast})^n) \subset \mathbb{P}(L(\Delta)^{\vee})$. Assume $0 \in Int(\Delta)$ and let $$ \Delta^{\circ} = \left \{ n \mid (m,n) \geq -1 \textrm{ for every } m \in \Delta \right \} \subset N_{\mathbb{R}}$$ be the \emph{polar} polytope of $\Delta$. A polytope $\Delta \subset M_{\mathbb{R}}$ with $ 0 \in Int(\Delta)$ is said to be \emph{reflexive} if $\Delta^{\circ} \subset N_{\mathbb{R}}$ is integral. A reflexive polytope $\Delta$ is said to be \emph{Fano} if every facet of $\Delta^{\circ}$ is the convex hall of a basis of $M$. Batyrev shows in \cite{Ba2} that $X_{\Delta}$ is a Fano variety if and only if $\Delta$ is reflexive and, in this case, the embedding
$i_{\Delta}$ is the anti-canonical embedding. The Fano variety $X_{\Delta}$ is smooth if and only if $\Delta^{\circ}$ is a Fano polytope.

\hspace{-0.6cm} Denote by $\Delta(k)$ the set of $k$-dimensional faces of $\Delta$ and let $ V_X(F) \subset X$ be the closure of the $T$-orbit corresponding to the facet $F \in \Delta(k)$ in $X$. In particular, consider the group of toric divisors $$ Div_T(X) := \bigoplus_{F \in \Delta(n-1)} \mathbb{Z} \cdot V_X(F).$$ When $X$ is a smooth toric
manifold the group $Pic(X)$ admits a description in terms of the following short exact sequence $$ 0 \rightarrow M \rightarrow
Div_T(X) \rightarrow Pic(X) \rightarrow 0. $$ When $X$ is also Fano one has $\Delta(n-1) \simeq \Delta^{\circ}(0)= \left \{n_F \vert F \in \Delta(n-1) \right \} \subset N_{\mathbb{R}}$ and the map on the left hand side is given by $ m \rightarrow \sum_F \left < m, n_F \right > \cdot V_X(F) $. We sometimes denote $V_X(n_F): = V_X(F)$. Moreover,
one can always assume that, up to a linear integral automorphism, $$\Delta^{\circ}(0) = \left \{ e_1,...,e_n, n_1,...,n_{\rho} \right \} \subset \mathbb{Z}^n. $$ In the considered cases one has $Pic(X) \simeq \bigoplus_{i=1}^{\rho} [V_X(n_i)] \cdot \mathbb{Z}$
and we write $$ [V_X(e_i)] := \sum_{j=1}^{\rho} e_{ij} \cdot [V_X(n_i)] \in Pic(X),$$ with $e_{ij} \in \mathbb{Z}$ for $1 \leq i \leq n$ and $1 \leq j \leq \rho$.

\hspace{-0.6cm} Let $\mathcal{D}^b(X)$ be the derived category of bounded complexes of
coherent sheaves on $X$, see \cite{GM,Tho}. An object $ E \in \mathcal{D}^b(X)$ is said to be \emph{exceptional} if $Hom(E,E)=\mathbb{C}$ and $Ext^i(E,E)=0$ for $i \neq 0$. An ordered collection $\left \{ E_1,...,E_N \right \} \subset \mathcal{D}^b(X)$
is said to be an \emph{exceptional collection} if each $E_j$ is exceptional and $$Ext^i(E_j,E_k) =0 \hspace{0.5cm} \textrm{ for } \hspace{0.25cm} j<k \textrm{ and all } i. $$ An
exceptional collection is said to be \emph{strongly exceptional} if also $Ext^i(E_j,E_k)=0$ for $j \leq k$ and $i \neq 0$. A strongly exceptional collection is called \emph{full} if its elements generate $\mathcal{D}^b(X)$ as a triangulated category. It should be noted, that results of Efimov show that there are toric Fano manifolds which do not admit any full strongly exceptional collection of line bundles, see \cite{E}.

\hspace{-0.6cm} Let $W \in L(\Delta^{\circ})$ be a Laurent polynomial whose Newton polytope is the Fano polytope $\Delta^{\circ}$. Let $Crit(W) \subset (\mathbb{C}^{\ast})^n$ be the solution scheme of the corresponding Landau-Ginzburg system of equations
$$ z_k \frac{ \partial}{ \partial z_k } W(z_1,...,z_n) = 0 \hspace{0.5cm} \textrm{ for } \hspace{0.25cm} k=1,...,n. $$ For a solution $z=(z^1,...,z^n) \in Crit(W)$ consider the $T$-divisor with real coefficents $$ D_W(z):=\sum_{F \in \Delta(n-1)} Arg(z^{n_F}) \cdot V_X(F) \in Div_T(X) \otimes \mathbb{R},$$ where $Arg : \mathbb{C}^{\ast} \rightarrow [0,1) $ is the Argument function given by $ Arg(r e^{2 \pi i \theta}) = \theta \in [0,1)$. Note that
$$ [D_W(z)]= \left [ \sum_{i=1}^n Arg(z^i) \cdot V_X(e_i) + \sum_{j=1}^{\rho} Arg(z^{n_i}) \cdot V_X(n_i) \right ] \in Pic(X) \otimes \mathbb{R}. $$ Define the map $E_W : Crit(W) \rightarrow Pic(X)$ as follows $$ E_W(z) := \left [ D_W(z) \right ]_{\mathbb{Z}} \in Pic(X) \hspace{0.5cm} \textrm{ for } \hspace{0.25cm} z \in Crit(W), $$ where $[D ]_{\mathbb{Z}} \in Pic(X)$ is the integral part of $[D] \in Pic(X) \otimes \mathbb{R}$. In \cite{J,J2} we showed various cases of pairs $(X,W)$ for which $E_W: Crit(W) \rightarrow Pic(X)$ is an \emph{exceptional} map, that is, its image is a full strongly exceptional collection. Let us briefly recall selected examples:
\bigskip

\hspace{-0.6cm} \bf Example 2.1 \rm (Projective space): \rm For $X = \mathbb{P}^n$ one has $\Delta^{\circ}(0) = \left \{ e_1,...,e_n, e_0 \right \}$ with $e_0:= - \sum_{i=1}^n e_i$ and $$ H := [V_X(e_0)]=[V_X(e_1)] =...= [V_X(e_n)] \in Pic(X).$$ The Landau-Ginzburg system of equations is given by $$ z_i \frac{ \partial }{\partial z_i} W(z_1,...,z_n) = z_i - \frac{1}{z_1 \cdot ... \cdot z_n} =0 \hspace{0.5cm} \textrm{ for } i=1,...,n,$$ whose solution scheme is given by $$ Crit(W) = \left \{ \left ( e^{\frac{2 \pi i k}{n+1}} ,..., e^{\frac{2 \pi i k}{n+1}} \right ) \vert k=0,...,n \right \} \subset (\mathbb{C}^{\ast})^n.$$ Hence, for $k=0,...,n$ one has $$E_W(z_k) =\left [ \sum_{i=0}^n \frac{k}{n+1} \cdot V_X(e_i) \right ]_{\mathbb{Z}} = \left [ \left (\sum_{i=0}^n \frac{k}{n+1} \right ) \cdot H \right ]_{\mathbb{Z}} = k \cdot H.$$ In particular, the resulting collection $\mathcal{E}_W=\left \{ k \cdot H \vert k = 0,...,n \right \} \subset Pic(X)$ is Beilinson's full strongly exceptional collection, see \cite{B}.

\bigskip

\hspace{-0.6cm} \bf Example 2.2 \rm (Products): \rm Let $X_1,X_2$ be two toric Fano manifolds given by polytopes $\Delta_1 \subset \mathbb{R}^n$ and $\Delta_2 \subset \mathbb{R}^m$. Set $$ \begin{array}{ccc} \Delta_1^{\circ}(0)=\left \{ e_1,...,e_n,n_1,...,n_{\rho_1} \right \} \subset \mathbb{Z}^n & ; & \Delta^{\circ}_2(0) = \left \{ v_1,...,v_m, n'_1,...,n'_{\rho_2} \right \} \subset \mathbb{Z}^m \end{array}. $$ Then $X:=X_1 \times X_2 $ is also a toric Fano manifold given by the polytope $ \Delta \subset \mathbb{R}^n \times \mathbb{R}^m$ satisfying $$ \Delta^{\circ}(0) = pr_1^{\ast} \Delta_1^{\circ}(0) \times pr^{\ast}_2 \Delta^{\circ}_2(0) \subset \mathbb{Z}^n \times \mathbb{Z}^m,$$ where $pr_1 : \mathbb{Z}^n \times \mathbb{Z}^m \rightarrow \mathbb{Z}^n$ and $ pr_2 : \mathbb{Z}^n \times \mathbb{Z}^m \rightarrow \mathbb{Z}^m$ are the two projections. In particular, one has $$ Pic(X) \simeq \pi_1^{\ast} Pic(X_1) \oplus \pi_2^{\ast} Pic(X_2).$$ If $\mathcal{E}_1 \subset Pic(X_1)$ and $\mathcal{E}_2 \subset Pic(X_2)$ are full strongly exceptional collections on $X_1$ and $X_2$, then $\mathcal{E}:=\pi_1^{\ast} \mathcal{E}_1 \oplus \pi_2^{\ast} \mathcal{E}_2 \subset Pic(X)$ is a full strongly exceptional collection on $X$, see \cite{CMR}. Assume $W_1(z) \in L(\Delta^{\circ}_1)$ and $W_2(w) \in L(\Delta^{\circ}_2)$ are such that $E_{W_i}: Crit(W_i) \rightarrow Pic(X_i)$ is an exceptional map for $i=1,2$. Set $W(z,w):=W_1(z)+W_2(w) \in L(\Delta^{\circ})$ and note that $$ Crit(W) = \left \{ (z,w) \vert
z \in Crit(W_1) , w \in Crit(W_2) \right \} \subset (\mathbb{C}^{\ast})^n \times (\mathbb{C}^{\ast})^m.$$ In particular, $E_W$ is an exceptional map on $X_1 \times X_2$ with $\mathcal{E}_W = \pi^{\ast}_1 E_{W_1} \oplus \pi^{\ast}_2 E_{W_2}$.
\bigskip

\hspace{-0.6cm} \bf Example 2.3 \rm (Toric Del Pezzo surfaces): In $dim(X)=2$ the only toric Fano manifolds are the five toric Del-Pezzo surfaces $$ X = \mathbb{P}^2,\mathbb{P}^1 \times \mathbb{P}^1, Bl_1(\mathbb{P}^2),Bl_2(\mathbb{P}^2),Bl_3(\mathbb{P}^2),$$ where $Bl_k(\mathbb{P}^2)$ is the blow-up of $\mathbb{P}^2$ in $k=1,2,3$ of the $T$-invariant points. The first two cases were given in Example 2.1 and Example 2.2. For the blow-ups, one has $$ Pic(Bl_k(\mathbb{P}^2)) \simeq H \cdot \mathbb{Z} \oplus \left ( \bigoplus_{i=1}^k E_i \cdot \mathbb{Z} \right ) ,$$ where $ E_i$ is the class of the normal bundle of the $i$-th exceptional divisor. It was shown in \cite{K,P} that $$ \mathcal{E} = \left \{ 0,H,H-E_i,2H-\sum_{i=1}^k E_i \vert i=1,...,k \right \} \subset Pic(Bl_k(\mathbb{P}^2))$$ is a full strongly exceptional collection on $Bl_k(\mathbb{P}^2)$ for $k=1,2,3$.

\bigskip

\hspace{-0.6cm} \underline{$k=1$:} The case $k=1$ is a special case of the following Example 2.4.

\bigskip

\hspace{-0.6cm} \underline{$k=2$:} For $k=2$ one has $ \Delta^{\circ}(0)= \left \{ e_1,e_2,n_1,n_2,n_3 \right \} $ with $$ \begin{array}{ccccc} n_1:=-e_1 & ; & n_2:=-e_1-e_2 & ; & n_3 = -e_2 \end{array}$$ with $$ \begin{array}{cccccc} [V_X(e_1)]=H-E_2 & ; & [V_X(e_2)]= H-E_1 & ; & [V_X(n_1)]=E_1 & ; \end{array} $$ $$ \begin{array}{ccc} [V_X(n_2)] = H-E_1-E_2 & ; & [V_X(n_3)]=E_2 \end{array}. $$ On the other hand, set $$ W(z_1,z_2) := z_1+z_2 + \frac{1}{z_1} + \frac{1}{z_1 z_2} + \frac{1}{z_2} \in L(\Delta^{\circ}).$$
In this case (as in the case $k=1$) computation shows that the Arguments of the solutions for the potential $W$ itself, are not given in terms of roots of unity. This led us in \cite{J2} to define for $0 \leq t$ the following deformation of the potential $$ W_t(z_1,z_2) :=e^{-t} z_1+e^{-t} z_2 + \frac{1}{z_1} + \frac{1}{z_1 z_2} + \frac{1}{z_2} \in L(\Delta^{\circ}).$$ Denote by $$\widetilde{z} := Arg (z^{e_1},z^{e_2},z^{n_1},z^{n_2},z^{n_3} ) \in \mathbb{T}^5.$$ Direct computation gives that $\widetilde{z}$ for the elements of the solution scheme $Crit(W_t) \subset (\mathbb{C}^{\ast})^2$ converge, as $ t \rightarrow \infty$, to the following five points $$ \begin{array}{cccccc} \widetilde{z}^0=(0,0,0,0,0) & ; & \widetilde{z}^1=(1/2,0,1/2,1/2,0) & ; & \widetilde{z}^2=(0,1/2,0,1/2,1/2) & ; \end{array} $$
$$ \begin{array}{ccc} \widetilde{z}^3=(1/2,1/2,1/2,0,1/2) & ; & \widetilde{z}^4=(1/2,1/2,1/2,1,1/2) \end{array}.$$ In particular, we have $E_W(\widetilde{z}^0)=0$ and $$ E_W(\widetilde{z}^1)= \left [ \frac{1}{2} [V_X(e_1)]+\frac{1}{2} [V_X(n_1)] + \frac{1}{2} [V_X(n_2)] \right]_{\mathbb{Z}}= $$ $$=\left [ \frac{1}{2} (H-E_2)+\frac{1}{2} E_1+ \frac{1}{2}(H-E_1-E_2) \right]_{\mathbb{Z}}=H-E_2$$

$$ E_W(\widetilde{z}^2)= \left [ \frac{1}{2} [V_X(e_2)]+\frac{1}{2} [V_X(n_2)] + \frac{1}{2} [V_X(n_3)] \right]_{\mathbb{Z}}= $$ $$=\left [ \frac{1}{2} (H-E_1)+\frac{1}{2} (H-E_1-E_2)+ \frac{1}{2}E_2 \right]_{\mathbb{Z}}=H-E_1$$

$$ E_W(\widetilde{z}^3)= \left [ \frac{1}{2} [V_X(e_1)]+\frac{1}{2} [V_X(e_2)]+ \frac{1}{2} [V_X(n_1)] + \frac{1}{2} [V_X(n_3)] \right]_{\mathbb{Z}}= $$
$$ \left [ \frac{1}{2} (H-E_2) + \frac{1}{2} (H-E_1) + \frac{1}{2} E_1 + \frac{1}{2} E_2 \right ]_{\mathbb{Z}} = H$$

$$ E_W(\widetilde{z}^4)= \left [ \frac{1}{2} [V_X(e_2)]+ \frac{1}{2} [V_X(e_2)] + \frac{1}{2} [V_X(n_1)] + [V_X(n_2)]+ \frac{1}{2}[V_X(n_3)] \right]_{\mathbb{Z}}=$$
$$ \left [\frac{1}{2} (H-E_2) + \frac{1}{2} (H-E_1) + \frac{1}{2} E_1 + (H-E_1-E_2) + \frac{1}{2} E_2 \right ]_{\mathbb{Z}} = 2H-E_1-E_2.$$

\bigskip

\hspace{-0.6cm} \underline{$k=3$:} For $k=3$ one has $ \Delta^{\circ}(0)= \left \{ e_1,e_2,n_1,n_2,n_3,n_4 \right \} $ with $$ \begin{array}{ccccccc} n_1:=e_1+e_2 & ; & n_2:=-e_1 & ; & n_3:=-e_1-e_2 & ; & n_4 = -e_2 \end{array}$$ with $$ \begin{array}{cccccc} [V_X(e_1)]=H-E_1-E_2 & ; & [V_X(e_2)]= H-E_1-E_2 & ; & [V_X(n_1)]=E_3 & ; \end{array} $$ $$ \begin{array}{ccccc} [V_X(n_2)] = E_1 & ; & [V_X(n_3)]=H-E_1-E_2 & ; & [V_X(n_4)] = E_2 \end{array}. $$ On the other hand, set $$ W(z_1,z_2) := z_1+z_2 + z_1 z_2 + \frac{1}{z_1} + \frac{1}{z_1 z_2} + \frac{1}{z_2} \in L(\Delta^{\circ}).$$ Direct computation gives that the elements of the solution scheme $Crit(W) \subset (\mathbb{C}^{\ast})^2$ are the following six points $$ \begin{array}{cccccccc} z^0=(1,1) & ; & z^1=(\rho_3,\rho_3) & ; & z^2=(\rho_3^2,\rho_3^2) & ; & z^3 =(1,-1) & ; \end{array} $$
$$ \begin{array}{ccc} z^4=(-1,1) & ; & z^5=(-1,-1) \end{array}$$ where $ \rho_3:=e^{\frac{2 \pi i}{3}} \in \mathbb{C}^{\ast}$. In particular, we have $E_W(z^0)=0$ and $$ E_W(z^1)= \left [ \frac{1}{3} [V_X(e_1)]+\frac{1}{3} [V_X(e_2)]+ \frac{2}{3} [V_X(n_1)] + \frac{2}{3} [V_X(n_2)]+ \frac{1}{3} [V_X(n_3)] + \frac{2}{3} [V_X(n_4)] \right]_{\mathbb{Z}}= $$ $$=\left [ \frac{1}{3} (H-E_2-E_3)+\frac{1}{3} (H-E_1-E_3)+ \frac{2}{3} E_3 + \frac{2}{3} E_1+ \frac{1}{3} (H-E_1-E_2) + \frac{2}{3} E_2 \right]_{\mathbb{Z}}=H$$

$$ E_W(z^2)= \left [ \frac{2}{3} [V_X(e_1)]+\frac{2}{3} [V_X(e_2)]+ \frac{1}{3} [V_X(n_1)] + \frac{1}{3} [V_X(n_2)]+ \frac{2}{3} [V_X(n_3)] + \frac{1}{3} [V_X(n_4)] \right]_{\mathbb{Z}}= $$ $$=\left [ \frac{2}{3} (H-E_2-E_3)+\frac{2}{3} (H-E_1-E_3)+ \frac{1}{3} E_3 + \frac{1}{3} E_1+ \frac{2}{3} (H-E_1-E_2) + \frac{2}{3} E_2 \right]_{\mathbb{Z}}=$$ $$=2H-E_1-E_2-E_3$$

$$ E_W(z^3)= \left [ \frac{1}{2} [V_X(e_2)]+ \frac{1}{2} [V_X(n_1)] + \frac{1}{2} [V_X(n_3)] + \frac{1}{2} [V_X(n_4)] \right]_{\mathbb{Z}}=$$ $$=\left [ \frac{1}{2} \left ( (H-E_1-E_3)+E_3+(H-E_1-E_2)+E_2 \right ) \right ]_{\mathbb{Z}}= H-E_1 $$

$$ E_W(z^4)= \left [ \frac{1}{2} [V_X(e_1)]+ \frac{1}{2} [V_X(n_1)] + \frac{1}{2} [V_X(n_2)] + \frac{1}{2} [V_X(n_3)] \right]_{\mathbb{Z}}=$$ $$=\left [ \frac{1}{2} \left ( (H-E_2-E_3)+E_3+E_1+(H-E_1-E_2) \right ) \right ]_{\mathbb{Z}}= H-E_2 $$

$$ E_W(z^5)= \left [ \frac{1}{2} [V_X(e_1)]+ \frac{1}{2} [V_X(n_2)] + \frac{1}{2} [V_X(n_2)] + \frac{1}{2} [V_X(n_4)] \right]_{\mathbb{Z}}=$$ $$=\left [ \frac{1}{2} \left ( (H-E_2-E_3)+(H-E_1-E_3)+E_1+E_2 \right ) \right ]_{\mathbb{Z}}= H-E_3. $$

\bigskip

\hspace{-0.6cm} \bf Example 2.4 \rm (Projective bundles): The projective bundle $$X= \mathbb{P}(\mathcal{O}_{\mathbb{P}^s} \oplus \mathcal{O}_{\mathbb{P}^s}(a_1) \oplus ... \oplus \mathcal{O}_{\mathbb{P}^s}(a_r)) \hspace{0.5cm} \textrm{ with } \hspace{0.25cm} a_1 \leq ... \leq a_r $$ is toric and Fano if also $\sum_{i=1}^r a_i \leq s$. The corresponding polytope $\Delta \subset \mathbb{R}^{r+s}$ satisfies $$ \Delta^{\circ}(0) = \left \{ e_1,...,e_s,v_1,...,v_s,e_0,v_0 \right \} \subset \mathbb{Z}^s \oplus \mathbb{Z}^r,$$ with $e_0 = \sum_{i=1}^r a_i v_i - \sum_{i=1}^s e_i$ and $v_0 = - \sum_{i=1}^r v_i$. Geometrically, $\xi= [V_X(v_0)] $ is the class of the tautological line bundle while $\pi^{\ast} H = [V_X(e_0)] $ is the class of the pull back of the generator $H \in Pic(\mathbb{P}^s)$ via the bundle map. One has $[V_X(v_i)]=\pi^{\ast} H $ for $i=1,...,s$ and $ [ V_X(e_i)] = \xi - a_i \cdot \pi^{\ast} H $ for $i=1,...,r$. Exceptional collections of objects on projective bundles were studied by Orlov in \cite{Or}. In \cite{CMR2} Costa and Miro-Roig show that, under the Fano condition, $$ \mathcal{E} = \left \{ k \cdot \pi^{\ast} H + l \cdot \xi \right \}_{k,l=0}^{s,r} \subset Pic(X)$$ is a full strongly exceptional collection. On the other hand consider $$ W_t(z,w):= \sum_{i=1}^s z_i + \sum_{i=1}^r w_i+ e^{-t} \cdot \frac{\prod_{i=1}^r w_i^{a_i}}{\prod_{i=1}^s z_i} + \frac{1}{\prod_{i=1}^r w_i} \in L(\Delta^{\circ}).$$ For $t \in \mathbb{R}$. In \cite{J2} we showed that $E_{W_t}: Crit(W_t) \rightarrow Pic(X)$ is an exceptional map for $0<<t$ big enough, and $\mathcal{E}_{W_t}$ coincides with the full strongly exceptional collection of Costa and Miro-Roig.

\bigskip

\hspace{-0.6cm} \bf Remark 2.5 \rm (Base-point potential): Note that in the above examples a specific base-point potential $W$ was introduced. In fact, as we are considering the integer-part in the definition of the maps $E_W$, the space $L(\Delta^{\circ})$ is actually divided into chambers away from the set of those elements for which $\vert E_W (Crit(W)) \vert < \chi(X)$, such that the image of the map $E_W$ is fixed for a given chamber. In particular, the potentials of Examples 2.1-2.4 could be viewed as representatives of their corresponding chambers in $L(\Delta^{\circ})$.

\bigskip

\hspace{-0.6cm} A classical fact in toric geometry is that $Div_T(X) \simeq \bigoplus_{F \in \Delta(n-1)} V_X(F) \cdot \mathbb{Z}$ admits a dual description as the group of certain piece-wise linear functions on $N_{\mathbb{R}}$. Let $\Sigma=\Sigma(\Delta)$ be the fan determined by the polyotope $\Delta$, see \cite{F,O}. When $X$ is Fano the vertices of the polar polytope $\Delta^{\circ}(0)$ are exactly the primitive integral generators of the $1$-dimensional rays $\rho \in \Sigma(1)$. We say that a function $\psi: N_{\mathbb{R}} \rightarrow \mathbb{R}$ is a $\Sigma$-\emph{support function} if it is continuous, linear on each maximal cone $\sigma \in \Sigma(n)$ and $\psi(n) \in \mathbb{Z}$ for any $n \in \Delta^{\circ}(0)$. Denote by $SF(N_{\mathbb{R}}, \Sigma)$ the group of all $\Sigma$-support functions. When $X$ is smooth one has $Div_T(X) \simeq SF(N_{\mathbb{R}}, \Sigma)$. Indeed, if $\psi \in SF(N_{\mathbb{R}},\Sigma)$ define $D_{\psi} := \sum_{n \in \Delta^{\circ}(0)} \psi(n) \cdot V_X(F_n)$. Given a $\Sigma$-support function $\psi \in SF(N_{\mathbb{R}}, \Sigma)$ and a maximal cone $\sigma \in \Sigma(n)$ denote by $m(\psi ; \sigma) \in M$ the unique element such that $\psi(n) = \left < m(\psi, \sigma) ,n \right> $ for any $n \in \sigma \cap N$. We similarly denote $m(D ; \sigma)=m(\psi_D ; \sigma)$ where $\psi_D$ is the support function such that $D= D_{\psi_D}$. In particular, in the Fano case for any $E \in Pic(X)$ and $ \sigma \in \Sigma(n)$ define the weight $$m_X(E ; \sigma) := m(D(E) ; \sigma) \in \mathbb{Z}^n,$$ where $$ D(E) := \sum_{i=1}^{\rho} \widetilde{a}(E ; \rho) \cdot V_X(n_i) \in Div_T(X)$$ is the unique divisor representing $E$ in the generators $[V_X(n_i)]$ for $i = 1,...,\rho$.

\section{Review of homological mirror symmetry for toric Fano manifolds}
\label{s:intro}

\hspace{-0.6cm} Our aim in this section is to review the relevant features of Abouzaid's mirror symmetry functor for toric manifolds, developed in \cite{Ab1,Ab2}. Let $W \in L(\Delta^{\circ})$ be a generic Laurent polynomial whose Newton polytope is $\Delta^{\circ}$. Denote by $M_W :=W^{-1}(0) \subset (\mathbb{C}^{\ast})^n$ the fiber of $W$ over $0 \in \mathbb{C}$. In \cite{Ab1} Abouzaid introduces the following family of functions $$\widetilde{W}_{t,s} (z) := 1 + \sum_{n \in \Delta^{\circ}(0)} t^{-1}(1-s \cdot \phi_n(Log \vert z \vert)) z^n,$$ where $\phi_n \in \mathbb{C}^{\infty}(\mathbb{R}^n)$ for $n \in \Delta^{\circ}(0)$ are required to satisfy certain decay conditions. For $s =0$ one has $\widetilde{W}_{t,0} \in L(\Delta^{\circ})$ but for $ s \neq 0$ the function $\widetilde{W}_{t,s}$ is no longer a Laurent polynomial. However, $M_{t,s}:= \widetilde{W}^{-1}(0)$ are all symplectomorphic for generic values of $s$ and $0<<t$ (see Proposition 4.9 of \cite{Ab1}) . Denote by $M=M_{t,1}$ with $0<<t$ big enough. We refer to the pair $((\mathbb{C}^{\ast})^n , M)$ as the \emph{(tropical) polarized mirror model} of the toric Fano manifold $X$.

\bigskip

\hspace{-0.6cm} \bf Definition 3.1: \rm A Lagrangian brane $L \subset (\mathbb{C}^{\ast})^n$ is an embedded compact graded Lagrangian submanifold, which is
spin and exact. A Lagrangian brane is said to be admissible if $ \partial L \subset M$ and there exists a small neighborhood of $\partial L$ in $L$ which agrees with the parallel transport of $\partial L$ along
a segment in $\mathbb{C}$. A pair of admissible Lagrangians $(L_1,L_2)$ is called positive if their corresponding segments $\gamma_1,\gamma_2 \subset \mathbb{C}$ lie in the left half-plane and their tangent
vectors are oriented counter clockwise such that $Im(\gamma_2(\theta))<Im(\gamma_1(\theta))$.

\bigskip

\hspace{-0.6cm} Admissible Lagrangian branes are objects of the Fukaya $\mathcal{A}_{\infty}$ pre-category $Fuk((\mathbb{C}^{\ast})^n,M)$ which we hence-forward denote $Fuk(X^{\circ})$, see \cite{Ab2,Ko2,S3}. In particular, the space of morphisms
between two positive transverse objects $L_1,L_2 \in Fuk(X^{\circ})$ is given by the Floer complex $ (CF^{\ast}(L_1,L_2), \partial)$. In \cite{Ab2} Abouzaid described the $ \mathcal{A}_{\infty}$ sub-pre-category of tropical Lagrangian sections
$$ Fuk_{trop}(X^{\circ}) \subset Fuk(X^{\circ}),$$ which he proved to be quasi-equivalent to the $DG$-category of line bundles on $X$.

\hspace{-0.6cm} Consider the map $Log: (\mathbb{C}^{\ast})^n \rightarrow \mathbb{R}^n$ and let $\mathcal{A} = \frac{1}{t} Log \vert M \vert \subset \mathbb{R}^n$ be the amoeba of $M$, see \cite{GKZ}. In \cite{Ab1,Ab2} Abouzaid shows that there exists a
component $\widetilde{\Delta} \subset \mathbb{R}^n \setminus \mathcal{A}$ in the complement of $\mathcal{A}$ which is contained in the polytope $\Delta \subset \mathbb{R}^n$ and is $C^0$-close to it. Note that the map $Log$ can be viewed as a fibration whose fiber is $\mathbb{T}^n$ and whose zero-section is the Lagrangian $(\mathbb{R}^+)^n \subset (\mathbb{C}^{\ast})^n$. Consider the following definition:

\bigskip

\hspace{-0.6cm} \bf Defintion 3.2: \rm A tropical Lagrangian section in $((\mathbb{C}^{\ast})^n,M)$ is an admissible Lagrangian brane $L$ which is a section of the map $Log$ restricted to $\widetilde{\Delta}$.

\bigskip

\hspace{-0.6cm} It is shown in \cite{Ab2} that up to Hamiltonian isotopy tropical Lagrangian sections in $((\mathbb{C}^{\ast})^n, M)$ are in one to one correspondence with elements of $Pic(X)$. The class of $\mathcal{O}_X \in Pic(X)$ corresponds to $L_0:=(\mathbb{R}^+)^n \cap Log^{-1}(\Delta)$, which is a tropical Lagrangian section of the polarized model.

\bigskip

\hspace{-0.6cm} \bf Example 3.3 \rm (the projective line): For $X = \mathbb{P}^1$ one has $\Delta = [-1,1] \subset \mathbb{R}$ and $Pic(X) \simeq \mathbb{Z}$. The mirror is given in this case by $(\mathbb{C}^{\ast},M_t )$ with
$M_t= \left\{e^{-t},e^t \right \} $. The inverse image of $\Delta$ under $Log_t$ is the annulus $e^{-t} \leq \vert z \vert \leq e^t$ in $\mathbb{C}^{\ast}$. In this case, one sees that tropical Lagrangian sections are curves in the annulus connecting the two points of $M$
passing through each circle $ \vert z \vert = const $ once. In particular, a section corresponds to the line bundle $\mathcal{O}(k) \in Pic(X)$ if it rotates around the origin $k \in \mathbb{Z}$ times.

\bigskip
\hspace{-0.6cm} Recall from section 2 that $Pic(X) \simeq SF(\Sigma)/\mathbb{Z}^n$ and that a $\Sigma$-support function $\psi \in SF(\Sigma)$ is determined by the elements $m(\psi ; \sigma) \in \mathbb{Z}^n$ for the maximal cones $\Sigma(n) \simeq \Delta(0)$. On the other hand, Abouzaid shows that any tropical Lagrangian section $L \subset (\mathbb{C}^{\ast})^n$ must coincide
with the zero section $L_0 \subset (\mathbb{C}^{\ast})^n$ in a small neighborhood of the fiber $Log^{-1}(x) \simeq \mathbb{T}^n$ for any vertex $x \in \Delta(0) \subset \mathbb{R}^n$. As $\mathbb{T}^n = \mathbb{R}^n/\mathbb{Z}^n$, a lift $\widetilde{L} \subset (\mathbb{R}^+)^n \times \mathbb{R}^n$ of the tropical section $L \subset (\mathbb{C}^{\ast})^n$ to the universal cover gives rise to a specification
of elements of $m(\widetilde{L} ; x) \in \mathbb{Z}^n$ to any $ x \in \Delta(0)$. In particular, one can define $\psi_{\widetilde{L}} \in SF(\Sigma)$ to be the $\Sigma$-support function given by $m(\psi ; \sigma) = m( \widetilde{L} ; x_{\sigma})$, for $x_{\sigma} \in \Delta(0) \simeq \Sigma(n)$. As a lift depends on a choice of a deck transformation $m \in \mathbb{Z}^n$ one has a well defined bijection $$ HMS : Pic(X) \rightarrow Fuk_{trop}(X^{\circ})/ \sim. $$ In \cite{Ab2} Abouzaid proves that this bijection extends, in fact, to an equivalence of $A_{\infty}$ pre-categories, where the pre-category on the left-hand side is the $DG$-category of line bundles on $X$. In particular, if $L_1,L_2 \in Fuk_{trop}(X^{\circ})$ are a positive pair, see \cite{Ab2}, then $$ HF(L_1,L_2) \simeq Hom(HMS(L_1),HMS(L_2))$$ as required.

\section{Landau-Ginzburg systems and tropical sections for toric Del-Pezzo}
\label{s:intro}

\hspace{-0.6cm} Let $X$ be a toric Fano manifold and let $W \in L(\Delta^{\circ})$ be the corresponding Landau-Ginzburg potential. Our aim is to describe, for various $X$, a map $$L : Crit(W) \rightarrow Fuk_{trop}(X^{\circ})$$ associating a tropical Lagrangian section $L(z) \subset (\mathbb{C}^{\ast})^n$ to a solution $z \in Crit(W)$ and satisfying $[L(z)] =HMS(E(z))$, where $HMS$ is Abouzaid's homological mirror symmetry functor. Let us start with the following observation: For $a=(a_1,...,a_n) \in \mathbb{R}^n$ consider $$ Crit_a(W) := \left \{ z_i \frac{ \partial }{\partial z_i} W(z_1,...,z_n) = a_i \vert i=1,...,n \right \} \subset (\mathbb{C}^{\ast})^n$$ and set $ Crit(W ; r ) : = \bigcup_{\vert a \vert \leq r} Crit_a(W) \subset (\mathbb{C}^{\ast})^n$ for $0 \leq r$. Geometrically, we think of $Crit(W ; r)$ as an inflation of the elements of $Crit(W) = Crit(W ; 0)$ into embedded Lagrangian $n$-balls of radius $r$. Concretely, for $z \in Crit(W)$ we have a Lagrangian embedding $i_z : B^n(r) \rightarrow (\mathbb{C}^{\ast})^n$ such that $i_z(0)=z$ and $Crit(W ; r) = \bigoplus_{z \in Crit(W)} \widetilde{L}_r(z)$ where $\widetilde{L}_r(z) := i_z(B^n(r))$.

\hspace{-0.6cm} In the considered examples, direct computation (or direct application of the definition of amobeas of algebraic hypersurfaces, see \cite{Mikh}) shows that $$ lim_{t \rightarrow \infty} \left (Log_t(\widetilde{L}_t(z)) \right ) = \Delta \subset \mathbb{R}^n$$ for any $z \in Crit(W)$. However, although they satisfy the condition of being sections of the $Log$ map over the polytope $\Delta$, the Lagrangians $\widetilde{L}_t(z)$ are not tropical sections as they are not, in general, admissible, that is, they do not satisfy $\partial \widetilde{L}_t(z) \subset M_t$. In fact, only $\widetilde{L}_t(z_0) \subset (\mathbb{C}^{\ast})^n$, where $z_0 \in Crit(W)$ is the solution for which $E(z_0) = \mathcal{O}_X$, turns to be a geniune tropical section.

\hspace{-0.6cm} Let $R_X \subset L(\Delta^{\circ})$ be the hypersurfaces of all $W$ such that $Crit(W)$ is non-reduced. Note that a loop $\gamma \subset L(\Delta^{\circ}) \setminus R_X$ gives rise to an autmorphism $ M(\gamma) : Crit(W) \rightarrow Crit(W) $ defined via analytic contionuation. We refer to the map $$ M : \pi_1( L(\Delta^{\circ}) \setminus R_X , W) \rightarrow Aut(Crit(W))$$ as the monodromy map, see \cite{J,J2}. The map $M$ naturally inflates to a map $$M : \pi_1( L(\Delta^{\circ}) \setminus R_X , W) \rightarrow Aut(Crit(W; r))$$ by similarly applying $M(\gamma): Crit_a(W) \rightarrow Crit_a(W)$ point-wise to the elements of $\widetilde{L}_t(z)$. For two solutions $z,z' \in Crit(W)$ let
$ G(z,z'):=\left \{ \gamma \vert M(\gamma)(z)=z' \right \}/ \sim $ be the set of all paths, up to homotopy, which take $z$ to $z'$. We have the following:

\bigskip

\hspace{-0.6cm} \bf Lemma 4.1: \rm Let $X$ be a toric Fano manifold and $W$ a potential such that $ W \not \in R_X$. Let $z \in Crit(W)$ be a solution and $\gamma \in G(z,z')$. Then $ M(\gamma) (\widetilde{L}_t(z)) = \widetilde{L}_t(z')$.

\bigskip

\hspace{-0.6cm} \bf Proof: \rm One has $Crit_a(W) = \left \{ z_1(a),...,z_N(a) \right \} $ where $z_i(a)$ is a continuous function of $a \in B^n(r)$. In particular, $\widetilde{L}_t(z_i)= \bigcup_{a \in B^n(r)} z_i(a) $ and $M_a(\gamma)$ change continuously as a function of $a$. If $M(\gamma) (z_i) = z_j$ then clearly $M_a(\gamma)(z_i(a))=z_j(a)$ for $a$ with $\vert a \vert $ small enough. Cover the $B^n(r)$ with small enough balls and proceed
accordingly. $\square$

\bigskip

\hspace{-0.6cm} Lemma 4.1 shows that fixed inflated monodromy $M(\gamma)$ of a path $\gamma \in G(z,z_0)$ always transfers the embedded Lagrangian ball $\widetilde{L}_r(z)$ to the trivial Lagrangian tropical section. We observe, however, that even though the resulting section is always $L_0$, different monodromies transfer $\widetilde{L}_r(z)$ to $L_0$ along different paths. In order to obtain a non-trivial tropical section we decompose $\widetilde{L}_r(z) \simeq \Delta$ into components, to each of whose boundary points we associate an element in $G(z,z_0)$. The resulting Lagrangian $L(z) $ would be defined by a collection of paths $ \gamma( z ; x_{\sigma}) \in G(z,z_0)$ associated to the vertices $x_{\sigma}$ of $\widetilde{L}_r(z)$, which are in turn determined by the polytope $\Delta$ itself. Let us start with:

\bigskip

\hspace{-0.6cm} \bf Example 4.2 \rm (The projective line): For $X= \mathbb{P}^1$ one has $ W(z):= z+\frac{1}{z} \in L(\Delta^{\circ})$. Direct computation shows $$ Crit(W ; r) =[-e^{s},-e^{-s}] \cup [e^{-s},e^s] \subset \mathbb{C}^{\ast}$$ where $r = 2 sinh(s)$. In particular, $Crit(W)=\left \{ z_0,z_1 \right \} \subset \mathbb{C}^{\ast}$ with $z_0=1$ and $z_1=-1$. Consider the two monodromies
$\gamma_+(\theta) := z+ \frac{e^{2 \pi i \theta}}{z}$ and $\gamma_-(\theta) := e^{2 \pi i \theta} z+ \frac{1}{z}$. Note the following illustration:

$$
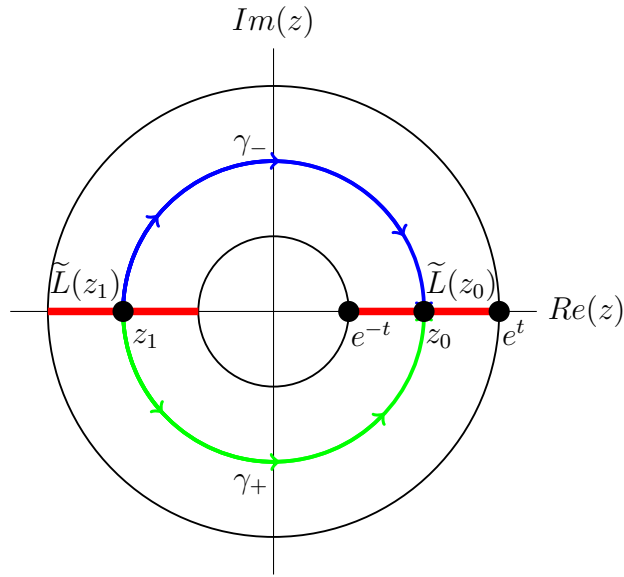
\begin{tikzpicture}

\draw (-3.5,0)--(3.5,0) node[right]{$Re(z)$};
\draw (0,-3.5)--(0,3.5) node[above]{$Im(z)$};
\draw [line width=0.25mm] (1,0) arc (0:360:1);

\draw (-2,0) [->,line width=0.5mm,blue] arc (180:87.5:2);
\draw (-2,0) [->,line width=0.5mm,blue] arc (180:140:2);
\draw (-2,0) [->,line width=0.5mm,blue] arc (180:30:2);
\draw (-2,0) [->,line width=0.5mm,blue] arc (180:0:2);

\draw (-2,0) [->,line width=0.5mm,green] arc (180:317.5:2);
\draw (-2,0) [->,line width=0.5mm,green] arc (180:222.5:2);
\draw (-2,0) [->,line width=0.5mm,green] arc (180:272.5:2);
\draw (-2,0) [->,line width=0.5mm,green] arc (180:360:2);

\draw [line width=0.25mm] (3,0) arc (0:360:3);

\draw [red,line width=1mm] (1,0) coordinate (a_1) -- (3,0) coordinate (a_2);
\draw [red,line width=1mm] (-1,0) coordinate (b_1) -- (-3,0) coordinate (b_2);
\node at (2,0) {$\bullet$};
\coordinate (a) at (2,0);
\coordinate (b) at (-2,0);

\coordinate (c) at (1,0);
\coordinate (d) at (3,0);

\coordinate (e) at (3.2,-0.28);
\coordinate (f) at (1.3,-0.3);

\coordinate (g) at (2.2,-0.36);
\coordinate (h) at (-1.7,-0.36);

\coordinate (i) at (-0.3,2.2);
\coordinate (j) at (-0.3,-2.3);

\coordinate (k) at (-2.5,0.4);
\coordinate (l) at (2.5,0.4);

\fill[black] (a) circle (4pt);
\fill[black] (b) circle (4pt);

\fill[black] (c) circle (4pt);
\fill[black] (d) circle (4pt);

\node at (e) {$e^t$};
\node at (f) {$e^{-t}$};

\node at (g) {$z_0$};
\node at (h) {$z_1$};

\node at (i) {$\gamma_-$};
\node at (j) {$\gamma_+$};

\node at (k) {$\widetilde{L}(z_1)$};
\node at (l) {$\widetilde{L}(z_0)$};
\end{tikzpicture}$$

\hspace{-0.6cm} The inscribed annulus is $Log^{-1}(\Delta) \subset \mathbb{C}^{\ast}$, as in example 3.3. The red segments are the two Lagrangians $\widetilde{L}(z_0), \widetilde{L}(z_1) \subset \mathbb{C}^{\ast}$. Note that the
Lagrangian $\widetilde{L}(z_0)$ coincides with the zero section. The blue and green lines describe the actions of the monodromies $M(\gamma_-),M(\gamma_+)$ respectively. Both monodromies take $\widetilde{L}(z_1)$ to $\widetilde{L}(z_0)$. However,
$M(\gamma_-)$ goes through the upper half-annulus while $M(\gamma_+)$ goes through the lower half-annulus. Going to the universal cover we get the following picture: $$
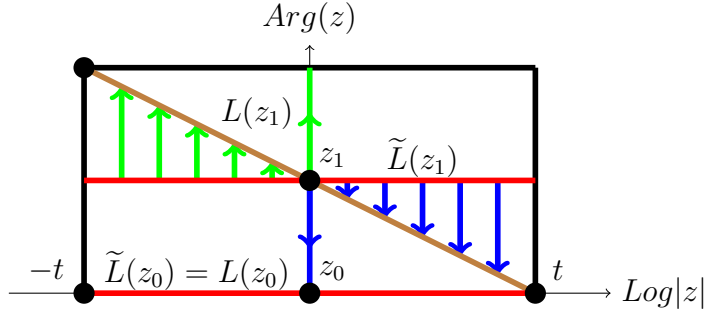
\begin{tikzpicture}
\draw [line width=0.75mm] (-3,3) coordinate (c_1) -- (3,3) coordinate (c_2);
\draw [line width=0.75mm] (-3,0) coordinate (d_1) -- (-3,3) coordinate (d_2);
\draw [line width=0.75mm] (3,0) coordinate (e_1) -- (3,3) coordinate (e_2);

\draw [->] (-4,0) coordinate (f_1) -- (4,0) coordinate (f_2) node[right]{$Log \vert z \vert$};
\draw [->] (0,0) coordinate (h_1) -- (0,3.3) coordinate (h_2) node[above] {$Arg(z)$};

\draw [line width=0.75mm,green] (0,1.5)--(0,3);
\draw [line width=0.75mm,blue] (0,1.5)-- (0,0);
\draw [->,line width=0.75mm,green] (0,1.5)--(0,2.4);
\draw [->,line width=0.75mm,blue] (0,1.5)-- (0,0.6);

\draw [->,line width=0.75mm,green] (-0.5,1.5)--(-0.5,1.75);
\draw [->,line width=0.75mm,green] (-1,1.5)--(-1,2);
\draw [->,line width=0.75mm,green] (-2,1.5)--(-2,2.5);
\draw [->,line width=0.75mm,green] (-1.5,1.5)--(-1.5,2.25);
\draw [->,line width=0.75mm,green] (-2.5,1.5)--(-2.5,2.75);

\draw [->,line width=0.75mm,blue] (0.5,1.5)--(0.5,1.25);
\draw [->,line width=0.75mm,blue] (1,1.5)--(1,1);
\draw [->,line width=0.75mm,blue] (2,1.5)--(2,0.5);
\draw [->,line width=0.75mm,blue] (1.5,1.5)--(1.5,0.75);
\draw [->,line width=0.75mm,blue] (2.5,1.5)--(2.5,0.25);

\draw [line width=0.75mm,brown] (3,0)--(-3,3);

\draw [red,line width=0.75mm] (-3,0) coordinate (a_1) -- (3,0) coordinate (a_2);
\draw [red,line width=0.75mm] (-3,1.5) coordinate (b_1) -- (3,1.5) coordinate (b_2);

\coordinate (b) at (0,1.5);

\fill[black] (h_1) circle (4pt);
\fill[black] (b) circle (4pt);
\fill[black] (d_1) circle (4pt);
\fill[black] (e_1) circle (4pt);
\fill[black] (-3,3) circle (4pt);

\node at (0.3,0.3) {$z_0$};
\node at (0.3,1.8) {$z_1$};

\node at (3.3,0.3) {$t$};
\node at (-3.5,0.3) {$-t$};

\node at (1.5,1.8) {$\widetilde{L}(z_1)$};
\node at (-1.5,0.3) {$\widetilde{L}(z_0)=L(z_0)$};
\node at (-0.7,2.4) {$L(z_1)$};

\end{tikzpicture}$$ We take $L(z_0)= \widetilde{L}(z_0)$ while $L(z_1)$ is the brown curve, obtained by applying $\gamma_-$ to the right part of $\widetilde{L}(z_1)$ and $\gamma_+$ to the left part of $\widetilde{L}(z_1)$.

\bigskip

\hspace{-0.6cm} Let $X$ be a toric Del Pezzo surface and let $ \sigma = \sum_{\rho < \sigma} \mathbb{R}^+ \cdot n_{ \rho} \in \Sigma(n)$ be a maximal cone. Denote by $$Ray^{\circ}(\sigma) = \left \{ \rho \vert \rho \not < \sigma \right \} \in \Sigma(1)$$ the set of all rays of $\Sigma$ which are not rays of $\sigma$. Note that $\vert Ray^{\circ}(\sigma) \vert = rank(Pic(X))$. In fact, it turns that $ \left \{ [V_X(n_{\rho})] \right \}_{\rho \in Ray^{\circ}(\sigma)} \subset Pic(X)$ is a set of generators for any $\sigma \in \Sigma(n)$. In particular, for any $E \in Pic(X)$ let $$ D(E ; \sigma) := \sum_{\rho \in Ray^{\circ}(\sigma) } a(E ; \rho ; \sigma) \cdot V_X(n_{\rho}) \in Div_T(X)$$ be the unique divisor such that $E=[D(E ; \sigma)] \in Pic(X)$ with $a(E, \rho ; \sigma)=0$ for $\rho < \sigma$. To a solution $z \in Crit(W)$ and a maximal cone $\sigma \in \Sigma(n)$ we associate the loop of Laurent polynomials $$ \gamma(z ; \sigma) (\theta) : = \sum_{\rho \in \Sigma(1)} e^{2 \pi i \theta \cdot a(E(z) ; \sigma ; \rho)} \cdot z^{n_{\rho}} \in L(\Delta^{\circ}).$$ Note that since $Arg(z_0)=0 \in \mathbb{T}^n \simeq \mathbb{R}^n / \mathbb{Z}^n$ an element $ \gamma \in G(z, z_0)$ can be associated with a weight $m( \gamma) \in \mathbb{Z}^n$ by applying $Arg : (\mathbb{C}^{\ast})^n \rightarrow \mathbb{T}^n$ to the corresponding path from $z$ to $z_0$ and
lifting to the universal cover. Set $m_W(z; \sigma) := m(\gamma(z ; \sigma)) \in \mathbb{Z}^n$. On the other hand, recall from the end of section 2 that an element $E \in Pic(X)$ is associated with the weight $m_X(E ; \sigma) \in \mathbb{Z}^n$. For toric Del-Pezzo surfaces we have:

\bigskip

\hspace{-0.6cm} \bf Proposition 4.3: \rm Let $X$ be a toric Del-Pezzo surface and let $W \in L(\Delta^{\circ})$ be the Landau-Ginzburg potential as in Example 2.3. Let $z \in Crit(W)$ be a solution. Then $\gamma(z ; \sigma) \in G(z,z_0)$ and $$m_W(z ;\sigma) = -m_X(E(z) ; \sigma)$$ for any maximal cone $\sigma \in \Sigma(2)$.

\bigskip

\hspace{-0.6cm} \bf Proof: \rm Our notations fo the solution shceme $Crit(W) = \left \{ z_0 , ..., z_N \right \} $, where $N= \chi(X)-1$, are as in Example 2.3. The proof requires the computation of the weights
$m_W(z ; \sigma)$ and $-m_X(E(z) ; \sigma)$ for any $\sigma \in \Sigma(2)$. The toric weights $m_X(E(z) ; \sigma)$ are computed via the standard formula. The monodromy weights $m_W(z ; \sigma)$ are computed directly. We include a few of the computations for $\mathbb{P}^2$ and $Bl_1(\mathbb{P}^2)$ to ilsutrate the general technique.
\bigskip

\hspace{-0.6cm} \underline{(1) Projective space:} The fan $\Sigma$ for $X=\mathbb{P}^2$ is given by
$$
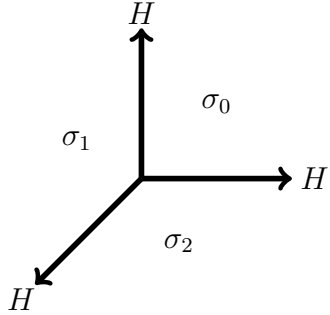
\begin{tikzpicture}

\draw [->,line width=0.75mm] (0,0) coordinate (a) -- (0,2) coordinate (b1);
\draw [->,line width=0.75mm] (0,0) coordinate (a) -- (2,0) coordinate (b2);
\draw [->,line width=0.75mm] (0,0) coordinate (a1) -- (-1.41,-1.41) coordinate (b3);

\node at (2.3,0) {$ H$};
\node at (0,2.2) {$H $};
\node at (-1.6,-1.6) {$H $};

\node at (1,1) {$\sigma_0$};
\node at (-0.86,0.5) {$\sigma_1$};
\node at (0.5,-0.86) {$\sigma_2$};

\end{tikzpicture}$$ Recall that in this case by $\mathcal{E}=\left \{ 0 , H,2H \right \} \subset Pic(\mathbb{P}^2) \simeq \mathbb{Z} \cdot H$. We describe the divisor $D(E) \in Div_T(X)$, in the sense of section 2, together with its corresponding weights as follows:

$$ \begin{array}{ccccc}

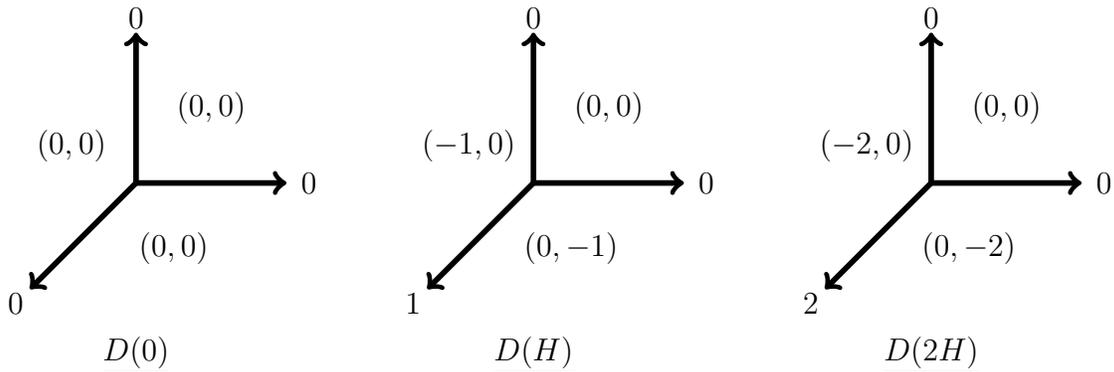
\begin{tikzpicture}

\draw [->,line width=0.75mm] (0,0) coordinate (a) -- (0,2) coordinate (b1);
\draw [->,line width=0.75mm] (0,0) coordinate (a) -- (2,0) coordinate (b2);
\draw [->,line width=0.75mm] (0,0) coordinate (a1) -- (-1.41,-1.41) coordinate (b3);

\node at (2.3,0) {$ 0 $};
\node at (0,2.2) {$0 $};
\node at (-1.6,-1.6) {$0 $};

\node at (1,1) {$(0,0)$};
\node at (-0.86,0.5) {$(0,0)$};
\node at (0.5,-0.86) {$(0,0)$};
\node at (0,-2.3) {$ \underline{D(0)}$};
\end{tikzpicture}
& \hspace{0.2cm} &
\begin{tikzpicture}

\draw [->,line width=0.75mm] (0,0) coordinate (a) -- (0,2) coordinate (b1);
\draw [->,line width=0.75mm] (0,0) coordinate (a) -- (2,0) coordinate (b2);
\draw [->,line width=0.75mm] (0,0) coordinate (a1) -- (-1.41,-1.41) coordinate (b3);

\node at (2.3,0) {$ 0 $};
\node at (0,2.2) {$0 $};
\node at (-1.6,-1.6) {$1 $};

\node at (1,1) {$(0,0)$};
\node at (-0.86,0.5) {$(-1,0)$};
\node at (0.5,-0.86) {$(0,-1)$};
\node at (0,-2.3) {$ \underline{D(H)}$};
\end{tikzpicture}

& \hspace{0.2cm} &

\begin{tikzpicture}

\draw [->,line width=0.75mm] (0,0) coordinate (a) -- (0,2) coordinate (b1);
\draw [->,line width=0.75mm] (0,0) coordinate (a) -- (2,0) coordinate (b2);
\draw [->,line width=0.75mm] (0,0) coordinate (a1) -- (-1.41,-1.41) coordinate (b3);

\node at (2.3,0) {$ 0$};
\node at (0,2.2) {$0 $};
\node at (-1.6,-1.6) {$2 $};
\node at (1,1) {$(0,0)$};
\node at (-0.86,0.5) {$(-2,0)$};
\node at (0.5,-0.86) {$(0,-2)$};
\node at (0,-2.3) {$ \underline{D(2H)}$};
\end{tikzpicture} \end{array}$$ On the other hand, let us consider the monodromies of $\gamma_i^j:=\gamma(z^i ; \sigma_j) \in G(z,z^0)$. We have $\gamma_0^j(\theta)=W$ for $j=1,..,3$, while
$$ \begin{array}{cccccc} \gamma_1^0(\theta) = z_1+z_2 + \frac{ e^{2 \pi i \theta }}{z_1 z_2} & ; & \gamma_1^1(\theta) = e^{2 \pi i \theta} z_1+z_2 + \frac{1}{z_1 z_2} & ;
& \gamma_1^2( \theta) = z_1+e^{2 \pi i \theta} z_2 + \frac{ 1}{z_1 z_2} \end{array} $$ and $$ \begin{array}{cccccc} \gamma_2^0(\theta) = z_1+z_2 + \frac{ e^{4 \pi i \theta }}{z_1 z_2} & ; & \gamma_2^1(\theta) = e^{4 \pi i \theta} z_1+z_2 + \frac{1}{z_1 z_2} & ; & \gamma_2^2(\theta) = z_1+e^{4 \pi i \theta} z_2 + \frac{ 1}{z_1 z_2} \end{array}. $$ Direct computation shows that the images of these paths, under the Argument map, lifted to the universal cover of $\mathbb{T}^2 = \mathbb{R}^2/\mathbb{Z}^2$ are given as follows:

$$
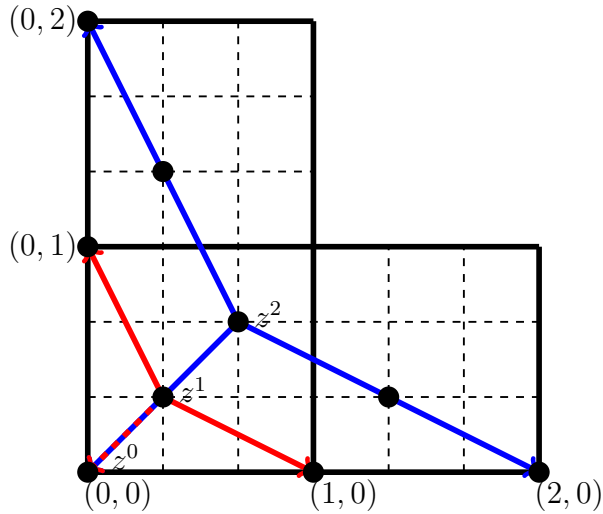
\begin{tikzpicture}

\draw [line width=0.75mm] (0,0) coordinate (a) -- (0,6) coordinate (b1);
\draw [line width=0.75mm] (0,0) coordinate (a) -- (6,0) coordinate (b2);
\draw [line width=0.75mm] (3,0) coordinate (a1) -- (3,6) coordinate (b3);
\draw [line width=0.75mm] (0,3) coordinate (a2) -- (6,3) coordinate (b4);
\draw [line width=0.75mm] (6,0) coordinate (b2) -- (6,3) coordinate (b4);
\draw [line width=0.75mm] (0,6) coordinate (b1) -- (3,6) coordinate (b3);

\draw [line width=0.25mm,dashed] (1,0)--(1,6);
\draw [line width=0.25mm,dashed] (2,0)-- (2,6);
\draw [line width=0.25mm,dashed] (0,1)-- (6,1);
\draw [line width=0.25mm,dashed] (0,2)--(6,2);
\draw [line width=0.25mm,dashed] (4,0)--(4,3);
\draw [line width=0.25mm,dashed] (5,0)-- (5,3);
\draw [line width=0.25mm,dashed] (0,4)--(3,4);
\draw [line width=0.25mm,dashed] (0,5)-- (3,5);

\draw [->,line width=0.75mm,blue] (2,2)--(0,0);
\draw [->,line width=0.75mm,blue] (2,2)--(6,0);
\draw [->,line width=0.75mm,blue] (2,2)--(0,6);

\draw [->,line width=0.75mm,dashed,red] (1,1)--(0,0);
\draw [->,line width=0.75mm,red] (1,1)--(3,0);
\draw [->,line width=0.75mm,red] (1,1)--(0,3);

\fill[black] (0,0) circle (4pt);
\fill[black] (1,1) circle (4pt);
\fill[black] (2,2) circle (4pt);
\fill[black] (0,3) circle (4pt);
\fill[black] (3,0) circle (4pt);
\fill[black] (0,6) circle (4pt);
\fill[black] (6,0) circle (4pt);
\fill[black] (4,1) circle (4pt);
\fill[black] (1,4) circle (4pt);

\node at (0.4,-0.3) {$ (0,0)$};
\node at (3.4,-0.3) {$(1,0)$};
\node at (6.4,-0.3) {$(2,0)$};
\node at (-0.6,3) {$(0,1)$};
\node at (-0.6,6) {$(0,2)$};
\node at (0.5,0.2) {$z^0$};
\node at (1.4,1.1) {$z^1$};
\node at (2.4,2.1) {$z^2$};

\end{tikzpicture}$$ As an illustration, let us compute the monodromy path corresponding to $\gamma^0_1(\theta)$ in order to do so we need to solve the system of equations $ z_1 - \frac{ e^{2 \pi i \theta}}{z_1 z_2}= z_2- \frac{e^{2 \pi i \theta}}{z_1 z_2} =0 $. Set $z=z_1=z_2$ and the second equation becomes $z^3= e^{2 \pi i \theta}$, which clearly agrees with the diagram.

\hspace{-0.6cm} Recall the $\mathbb{R}$-divisors $D_W(z) \in Div_T(X)$ defined in section 2. It is interesting to note that $m(D_W(z),\sigma) \in \mathbb{Z}^2$ is exactly the slope of $\gamma(z ; \sigma)$ under the Argument map. Indeed:
$$ \begin{array}{ccccc}

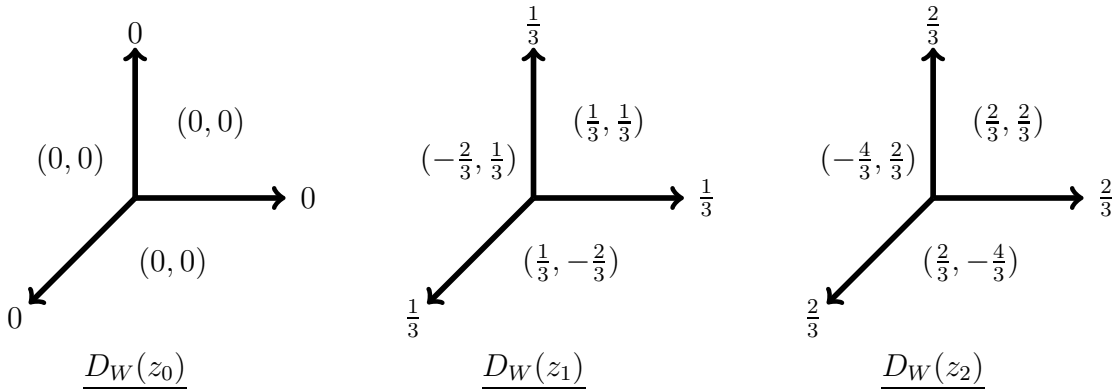
\begin{tikzpicture}

\draw [->,line width=0.75mm] (0,0) coordinate (a) -- (0,2) coordinate (b1);
\draw [->,line width=0.75mm] (0,0) coordinate (a) -- (2,0) coordinate (b2);
\draw [->,line width=0.75mm] (0,0) coordinate (a1) -- (-1.41,-1.41) coordinate (b3);

\node at (2.3,0) {$ 0 $};
\node at (0,2.2) {$0 $};
\node at (-1.6,-1.6) {$0 $};

\node at (1,1) {$(0,0)$};
\node at (-0.86,0.5) {$(0,0)$};
\node at (0.5,-0.86) {$(0,0)$};
\node at (0,-2.3) {$ \underline{D_W(z_0)}$};
\end{tikzpicture}
& \hspace{0.2cm} &
\begin{tikzpicture}

\draw [->,line width=0.75mm] (0,0) coordinate (a) -- (0,2) coordinate (b1);
\draw [->,line width=0.75mm] (0,0) coordinate (a) -- (2,0) coordinate (b2);
\draw [->,line width=0.75mm] (0,0) coordinate (a1) -- (-1.41,-1.41) coordinate (b3);

\node at (2.3,0) {$ \frac{1}{3} $};
\node at (0,2.3) {$\frac{1}{3} $};
\node at (-1.6,-1.6) {$\frac{1}{3} $};

\node at (1,1) {$(\frac{1}{3},\frac{1}{3})$};
\node at (-0.86,0.5) {$(-\frac{2}{3},\frac{1}{3})$};
\node at (0.5,-0.86) {$(\frac{1}{3},-\frac{2}{3})$};
\node at (0,-2.3) {$ \underline{D_W(z_1)}$};
\end{tikzpicture}

& \hspace{0.2cm} &

\begin{tikzpicture}

\draw [->,line width=0.75mm] (0,0) coordinate (a) -- (0,2) coordinate (b1);
\draw [->,line width=0.75mm] (0,0) coordinate (a) -- (2,0) coordinate (b2);
\draw [->,line width=0.75mm] (0,0) coordinate (a1) -- (-1.41,-1.41) coordinate (b3);

\node at (2.3,0) {$ \frac{2}{3}$};
\node at (0,2.3) {$\frac{2}{3} $};
\node at (-1.6,-1.6) {$\frac{2}{3} $};
\node at (1,1) {$(\frac{2}{3},\frac{2}{3})$};
\node at (-0.86,0.5) {$(-\frac{4}{3},\frac{2}{3})$};
\node at (0.5,-0.86) {$(\frac{2}{3},-\frac{4}{3})$};
\node at (0,-2.3) {$ \underline{D_W(z_2)}$};
\end{tikzpicture} \end{array}$$

\hspace{-0.6cm} \underline{(2) Blow up of $\mathbb{P}^2$ at one point:} The fan $\Sigma$ for $X=Bl_1(\mathbb{P}^2)$ is given by
$$
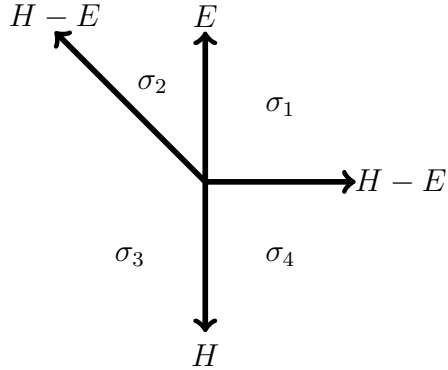
\begin{tikzpicture}

\draw [->,line width=0.75mm] (0,0) coordinate (a) -- (0,2) coordinate (b1);
\draw [->,line width=0.75mm] (0,0) coordinate (a) -- (2,0) coordinate (b2);
\draw [->,line width=0.75mm] (0,0)--(0,-2);
\draw [->,line width=0.75mm] (0,0) coordinate (a1) -- (-2,2) coordinate (b3);

\node at (2.6,0) {$ H-E $};
\node at (0,2.2) {$E $};
\node at (-2,2.2) {$H-E $};
\node at (0,-2.3) {$H $};

\node at (1,1) {$\sigma_1$};
\node at (1,-1) {$\sigma_4$};

\node at (-0.7,1.3) {$\sigma_2$};
\node at (-1,-1) {$\sigma_3$};

\end{tikzpicture}$$ We have $\mathcal{E}=\left \{ 0 , H,H-E,2H-E \right \} \subset Pic(Bl_1(\mathbb{P}^2)) \simeq \mathbb{Z} \cdot H \oplus \mathbb{Z} \cdot E $. The divisors $D(E) \in Div_T(X)$, together with the corresponding weights $m_X(E, \sigma)$, are given as follows:

$$ \begin{array}{ccc}
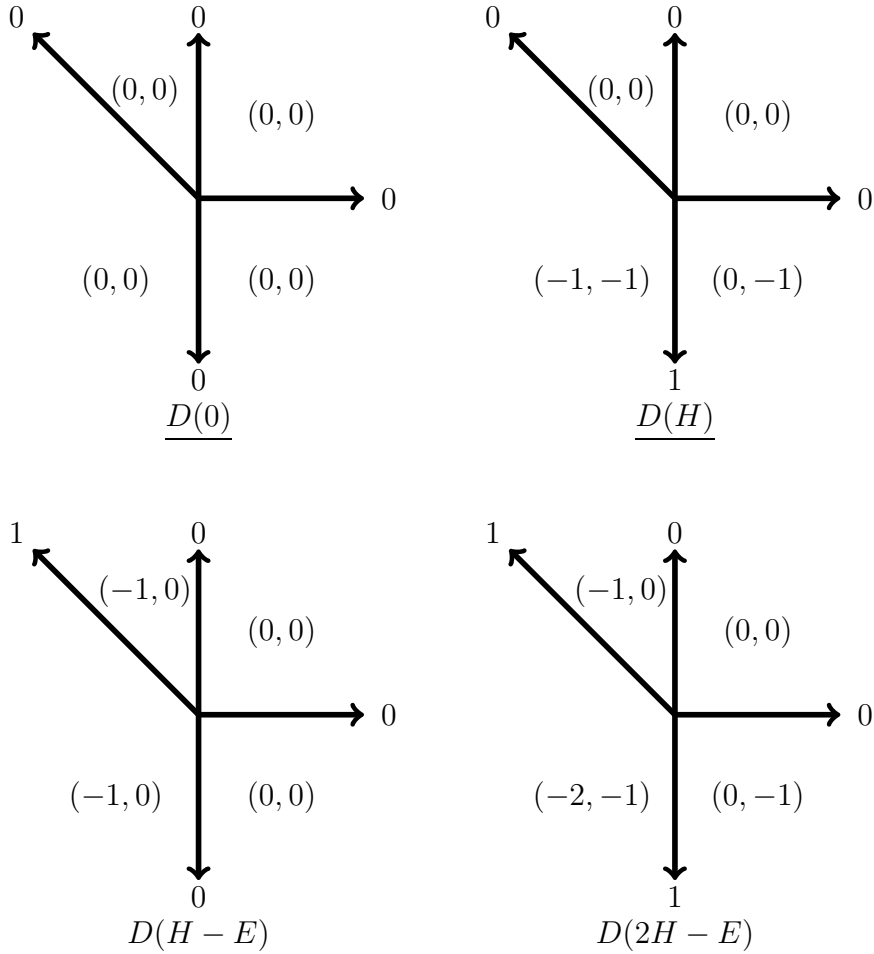
\begin{tikzpicture} [scale=1.1]

\draw [->,line width=0.75mm] (0,0) coordinate (a) -- (0,2) coordinate (b1);
\draw [->,line width=0.75mm] (0,0) coordinate (a) -- (2,0) coordinate (b2);
\draw [->,line width=0.75mm] (0,0)--(0,-2);
\draw [->,line width=0.75mm] (0,0) coordinate (a1) -- (-2,2) coordinate (b3);

\node at (2.3,0) {$ 0 $};
\node at (0,2.2) {$0 $};
\node at (-2.2,2.2) {$0 $};
\node at (0,-2.2) {$0 $};

\node at (1,1) {$(0,0)$};
\node at (1,-1) {$(0,0)$};

\node at (-0.65,1.3) {$(0,0)$};
\node at (-1,-1) {$(0,0)$};
\node at (0,-2.7) {$\underline{D(0)}$};
\end{tikzpicture}
& \hspace{0.2cm} &
\begin{tikzpicture}[scale=1.1]

\draw [->,line width=0.75mm] (0,0) coordinate (a) -- (0,2) coordinate (b1);
\draw [->,line width=0.75mm] (0,0) coordinate (a) -- (2,0) coordinate (b2);
\draw [->,line width=0.75mm] (0,0)--(0,-2);
\draw [->,line width=0.75mm] (0,0) coordinate (a1) -- (-2,2) coordinate (b3);

\node at (2.3,0) {$ 0 $};
\node at (0,2.2) {$0 $};
\node at (-2.2,2.2) {$0 $};
\node at (0,-2.2) {$1 $};

\node at (1,1) {$(0,0)$};
\node at (1,-1) {$(0,-1)$};

\node at (-0.65,1.3) {$(0,0)$};
\node at (-1,-1) {$(-1,-1)$};
\node at (0,-2.7) {$\underline{D(H)}$};

\end{tikzpicture} \end{array}$$
$$ \begin{array}{ccc}
\begin{tikzpicture}[scale=1.1]

\draw [->,line width=0.75mm] (0,0) coordinate (a) -- (0,2) coordinate (b1);
\draw [->,line width=0.75mm] (0,0) coordinate (a) -- (2,0) coordinate (b2);
\draw [->,line width=0.75mm] (0,0)--(0,-2);
\draw [->,line width=0.75mm] (0,0) coordinate (a1) -- (-2,2) coordinate (b3);

\node at (2.3,0) {$ 0 $};
\node at (0,2.2) {$0 $};
\node at (-2.2,2.2) {$1 $};
\node at (0,-2.2) {$0 $};

\node at (1,1) {$(0,0)$};
\node at (1,-1) {$(0,0)$};

\node at (-0.65,1.5) {$(-1,0)$};
\node at (-1,-1) {$(-1,0)$};
\node at (0,-2.7) {$\underline{D(H-E)}$};
\end{tikzpicture}
& \hspace{0.2cm} &
\begin{tikzpicture}[scale=1.1]

\draw [->,line width=0.75mm] (0,0) coordinate (a) -- (0,2) coordinate (b1);
\draw [->,line width=0.75mm] (0,0) coordinate (a) -- (2,0) coordinate (b2);
\draw [->,line width=0.75mm] (0,0)--(0,-2);
\draw [->,line width=0.75mm] (0,0) coordinate (a1) -- (-2,2) coordinate (b3);

\node at (2.3,0) {$ 0 $};
\node at (0,2.2) {$0 $};
\node at (-2.2,2.2) {$1 $};
\node at (0,-2.2) {$1 $};

\node at (1,1) {$(0,0)$};
\node at (1,-1) {$(0,-1)$};

\node at (-0.65,1.5) {$(-1,0)$};
\node at (-1,-1) {$(-2,-1)$};
\node at (0,-2.7) {$\underline{D(2H-E)}$};

\end{tikzpicture} \end{array}$$ On the other hand, the monodromies $\gamma_{ij}^k:=\gamma(z_{ij} ; \sigma_k) \in G(z_{ij},z_{00})$ are given as follows: $\gamma_{00}^k(\theta)=W$ for $k=1,...,4$ while
$$ \begin{array}{ccc} \gamma_{10}^1(\theta)= z_1+z_2+e^{-t} \cdot \frac{z_2}{z_1}+\frac{e^{2 \pi i \theta}}{z_2} & ; & \gamma_{10}^2(\theta)= z_1+z_2+e^{-t} \cdot \frac{z_2}{z_1}+\frac{e^{2 \pi i \theta}}{z_2} \end{array}$$
$$ \begin{array}{ccc} \gamma_{10}^3(\theta)= e^{2 \pi i \theta} \cdot z_1+e^{2 \pi i \theta} \cdot z_2+e^{-t} \cdot \frac{z_2}{z_1}+\frac{1}{z_2} & ; & \gamma_{10}^2(\theta)= z_1+e^{2 \pi i \theta } \cdot z_2+e^{-t+ 2 \pi i \theta} \cdot \frac{z_2}{z_1}+\frac{1}{z_2} \end{array}$$
$$ \begin{array}{ccc} \gamma_{01}^1(\theta)= z_1+z_2+e^{-t +2 \pi i \theta} \cdot \frac{z_2}{z_1}+ \frac{1}{z_2} & ; & \gamma_{01}^2(\theta)= e^{2 \pi i \theta} z_1+z_2+e^{-t} \cdot \frac{z_2}{z_1}+ \frac{1}{z_2} \end{array}$$
$$ \begin{array}{ccc} \gamma_{01}^3(\theta)= e^{2 \pi i \theta} \cdot z_1+z_2+e^{-t} \cdot \frac{z_2}{z_1}+ \frac{1}{z_2} & ; & \gamma_{01}^4(\theta)= z_1+z_2+e^{-t+2 \pi i \theta} \cdot \frac{z_2}{z_1}+ \frac{1}{z_2} \end{array}$$
$$ \begin{array}{ccc} \gamma_{11}^1(\theta)= z_1+z_2+e^{-t+2 \pi i \theta} \cdot \frac{z_2}{z_1}+e^{2 \pi i \theta} \frac{1}{z_2} & ; & \gamma_{11}^2(\theta)= e^{2 \pi i \theta} \cdot z_1+z_2+e^{-t} \cdot \frac{z_2}{z_1}+\frac{e^{2 \pi i \theta}}{z_2} \end{array}$$
$$ \begin{array}{ccc} \gamma_{11}^3(\theta)= e^{4 \pi i \theta} \cdot z_1+e^{ 2 \pi i \theta} \cdot z_2+e^{-t} \cdot \frac{z_2}{z_1}+\frac{1}{z_2} & ; & \gamma_{11}^4(\theta)= z_1+e^{2 \pi i \theta} \cdot z_2+e^{-t+4 \pi i \theta} \cdot \frac{z_2}{z_1}+e^{2 \pi i \theta} \frac{1}{z_2} \end{array}.$$ The images of these paths under the Argument map, lifted to the universal cover of $\mathbb{T}^2 = \mathbb{R}^2/\mathbb{Z}^2$ are:

$$
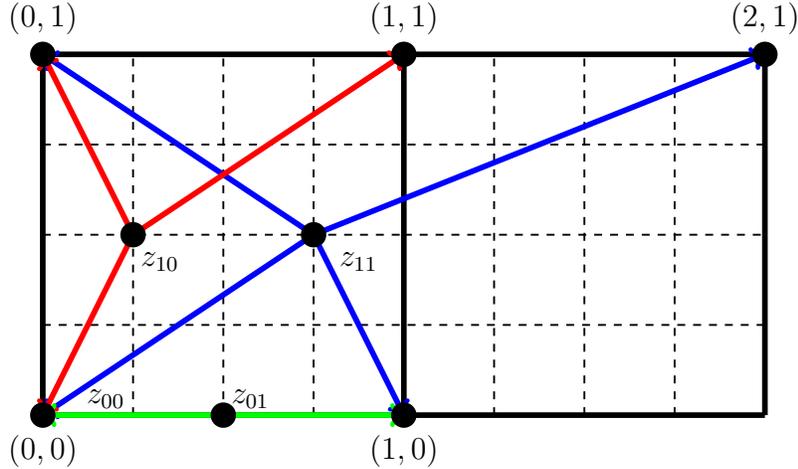
\begin{tikzpicture}[scale=1.2]

\draw [line width=0.75mm] (0,0)--(0,4);
\draw [line width=0.75mm] (0,0)--(8,0);
\draw [line width=0.75mm] (8,0)-- (8,4);
\draw [line width=0.75mm] (0,4)-- (8,4);
\draw [line width=0.75mm] (4,0)--(4,4);

\draw [line width=0.25mm,dashed] (2,0)--(2,4);
\draw [line width=0.25mm,dashed] (0,2)-- (4,2);
\draw [line width=0.25mm,dashed] (6,0)-- (6,4);
\draw [line width=0.25mm,dashed] (4,2)--(8,2);

\draw [line width=0.25mm,dashed] (1,0)--(1,4);
\draw [line width=0.25mm,dashed] (3,0)-- (3,4);
\draw [line width=0.25mm,dashed] (5,0)-- (5,4);
\draw [line width=0.25mm,dashed] (7,0)--(7,4);
\draw [line width=0.25mm,dashed] (0,1)--(8,1);
\draw [line width=0.25mm,dashed] (0,3)-- (8,3);

\draw [->,line width=0.75mm,blue] (3,2)--(8,4);
\draw [->,line width=0.75mm,blue] (3,2)--(0,4);
\draw [->,line width=0.75mm,blue] (3,2)--(0,0);
\draw [->,line width=0.75mm,blue] (3,2)--(4,0);

\draw [->,line width=0.75mm,red] (1,2)--(4,4);
\draw [->,line width=0.75mm,red] (1,2)--(0,4);
\draw [->,line width=0.75mm,red] (1,2)--(0,0);

\draw [->,line width=0.75mm,green] (2,0)--(4,0);
\draw [->,line width=0.75mm,green] (2,0)--(0,0);

\fill[black] (0,0) circle (4pt);
\fill[black] (1,2) circle (4pt);
\fill[black] (2,0) circle (4pt);
\fill[black] (3,2) circle (4pt);
\fill[black] (0,4) circle (4pt);
\fill[black] (4,0) circle (4pt);
\fill[black] (4,4) circle (4pt);
\fill[black] (8,4) circle (4pt);

\node at (2.3,0.18) {$z_{01}$};
\node at (1.3,1.7) {$z_{10}$};
\node at (3.5,1.7) {$z_{11}$};
\node at (0.7,0.18) {$z_{00}$};
\node at (0,-0.4) {$(0,0)$};
\node at (4,-0.4) {$(1,0)$};
\node at (0,4.4) {$(0,1)$};
\node at (4,4.4) {$(1,1)$};
\node at (8,4.4) {$(2,1)$};
\end{tikzpicture}$$ Let us compute the path corresponding to $\gamma^1_{10}(\theta)$. The corresponding system of equations is $ z_1 - e^{-t}\frac{z_2}{ z_1}= z_2+ e^{-t}\frac{z_2}{ z_1}-\frac{e^{2 \pi i \theta}}{z_2}=0 $. Due to the degeneration as $ t \rightarrow \infty$ this becomes equivalent to $ z_1 - e^{-t}\frac{z_2}{ z_1}= z_2-\frac{e^{2 \pi i \theta}}{z_2}=0 $. Hence $z^2_2 = e^{2 \pi i \theta}$ and $z_1^2 = e^{-t} z_2$. In particular, the path lies on the line $2 \theta_1 = \theta_2$ which agrees with the diagram. The rest of the cases are computed similarly.

\bigskip

\hspace{-0.6cm} \underline{(3) Blow up of $\mathbb{P}^2$ at two points:} The fan $\Sigma$ for $X=Bl_2(\mathbb{P}^2)$ is given by
$$
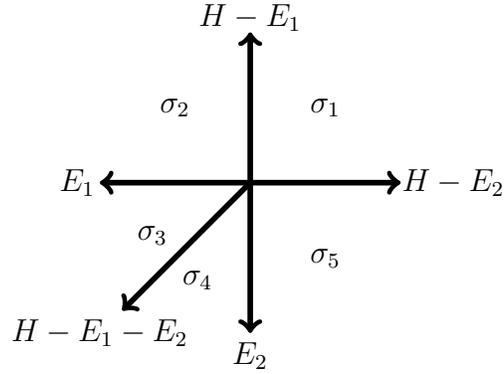
\begin{tikzpicture}

\draw [->,line width=0.75mm] (0,0) coordinate (a) -- (0,2) coordinate (b1);
\draw [->,line width=0.75mm] (0,0) coordinate (a) -- (2,0) coordinate (b2);
\draw [->,line width=0.75mm] (0,0) coordinate (a1) -- (-1.7,-1.7) coordinate (b3);
\draw [->,line width=0.75mm] (0,0)-- (0,-2);
\draw [->,line width=0.75mm] (0,0)-- (-2,0);

\node at (2.7,0) {$ H-E_2 $};
\node at (-2.3,0) {$E_1 $};
\node at (0,-2.3) {$E_2 $};
\node at (-2,-2) {$H-E_1-E_2 $};
\node at (0,2.2) {$H-E_1 $};

\node at (1,-1) {$\sigma_5$};
\node at (1,1) {$\sigma_1$};
\node at (-1.3,-0.7) {$\sigma_3$};

\node at (-0.7,-1.3) {$\sigma_4$};
\node at (-1,1) {$\sigma_2$};

\end{tikzpicture}$$ We have $\mathcal{E}=\left \{ 0 , H,H-E_1,H-E_2,2H-E_1-E_2 \right \} \subset Pic(Bl_2(\mathbb{P}^2))$. The divisors $D(E) \in Div_T(X)$, together with the corresponding weights $m_X(E, \sigma)$, are given as follows:
$$ 
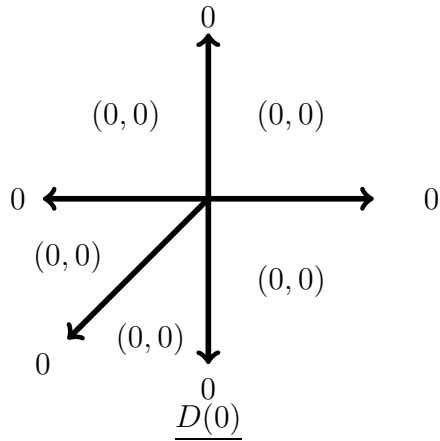
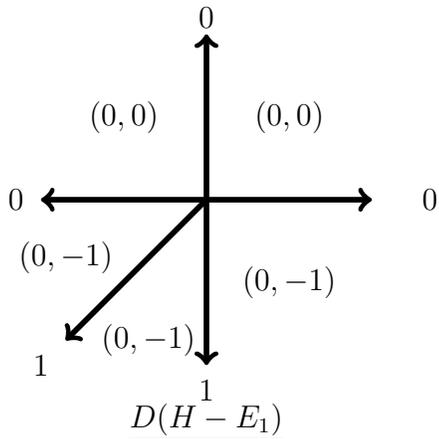
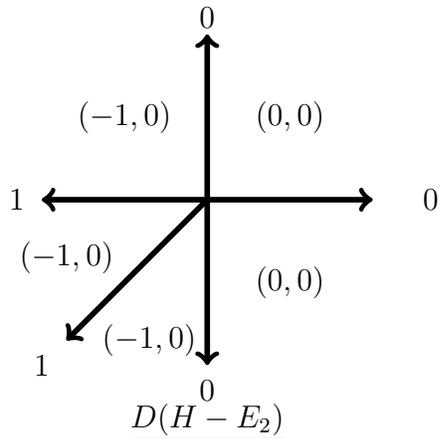
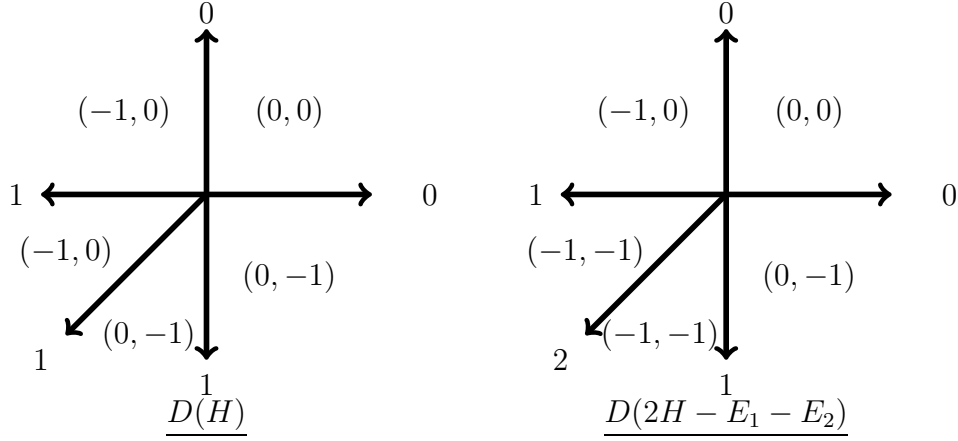
\begin{tikzpicture}[scale=1.1]

\draw [->,line width=0.75mm] (0,0) coordinate (a) -- (0,2) coordinate (b1);
\draw [->,line width=0.75mm] (0,0) coordinate (a) -- (2,0) coordinate (b2);
\draw [->,line width=0.75mm] (0,0) coordinate (a1) -- (-1.7,-1.7) coordinate (b3);
\draw [->,line width=0.75mm] (0,0)-- (0,-2);
\draw [->,line width=0.75mm] (0,0)-- (-2,0);

\node at (2.7,0) {$ 0 $};
\node at (-2.3,0) {$0 $};
\node at (0,-2.3) {$0 $};
\node at (-2,-2) {$0 $};
\node at (0,2.2) {$0 $};

\node at (1,-1) {$(0,0)$};
\node at (1,1) {$(0,0)$};
\node at (-1.7,-0.7) {$(0,0)$};

\node at (-0.7,-1.7) {$(0,0)$};
\node at (-1,1) {$(0,0)$};
\node at (0,-2.7) {$\underline{D(0)}$};
\end{tikzpicture}$$

$$ \begin{array}{ccc}
\begin{tikzpicture}[scale=1.1]

\draw [->,line width=0.75mm] (0,0) coordinate (a) -- (0,2) coordinate (b1);
\draw [->,line width=0.75mm] (0,0) coordinate (a) -- (2,0) coordinate (b2);
\draw [->,line width=0.75mm] (0,0) coordinate (a1) -- (-1.7,-1.7) coordinate (b3);
\draw [->,line width=0.75mm] (0,0)-- (0,-2);
\draw [->,line width=0.75mm] (0,0)-- (-2,0);

\node at (2.7,0) {$ 0 $};
\node at (-2.3,0) {$0 $};
\node at (0,-2.3) {$1 $};
\node at (-2,-2) {$1 $};
\node at (0,2.2) {$0 $};

\node at (1,-1) {$(0,-1)$};
\node at (1,1) {$(0,0)$};
\node at (-1.7,-0.7) {$(0,-1)$};

\node at (-0.7,-1.7) {$(0,-1)$};
\node at (-1,1) {$(0,0)$};
\node at (0,-2.7) {$\underline{D(H-E_1)}$};

\end{tikzpicture}
& \hspace{0.2cm} &
\begin{tikzpicture}[scale=1.1]

\draw [->,line width=0.75mm] (0,0) coordinate (a) -- (0,2) coordinate (b1);
\draw [->,line width=0.75mm] (0,0) coordinate (a) -- (2,0) coordinate (b2);
\draw [->,line width=0.75mm] (0,0) coordinate (a1) -- (-1.7,-1.7) coordinate (b3);
\draw [->,line width=0.75mm] (0,0)-- (0,-2);
\draw [->,line width=0.75mm] (0,0)-- (-2,0);

\node at (2.7,0) {$ 0 $};
\node at (-2.3,0) {$1 $};
\node at (0,-2.3) {$0 $};
\node at (-2,-2) {$1 $};
\node at (0,2.2) {$0 $};

\node at (1,-1) {$(0,0)$};
\node at (1,1) {$(0,0)$};
\node at (-1.7,-0.7) {$(-1,0)$};

\node at (-0.7,-1.7) {$(-1,0)$};
\node at (-1,1) {$(-1,0)$};
\node at (0,-2.7) {$\underline{D(H-E_2)}$};

\end{tikzpicture} \end{array}$$

$$ \begin{array}{ccc}
\begin{tikzpicture}[scale=1.1]

\draw [->,line width=0.75mm] (0,0) coordinate (a) -- (0,2) coordinate (b1);
\draw [->,line width=0.75mm] (0,0) coordinate (a) -- (2,0) coordinate (b2);
\draw [->,line width=0.75mm] (0,0) coordinate (a1) -- (-1.7,-1.7) coordinate (b3);
\draw [->,line width=0.75mm] (0,0)-- (0,-2);
\draw [->,line width=0.75mm] (0,0)-- (-2,0);

\node at (2.7,0) {$ 0 $};
\node at (-2.3,0) {$1 $};
\node at (0,-2.3) {$1 $};
\node at (-2,-2) {$1 $};
\node at (0,2.2) {$0 $};

\node at (1,-1) {$(0,-1)$};
\node at (1,1) {$(0,0)$};
\node at (-1.7,-0.7) {$(-1,0)$};

\node at (-0.7,-1.7) {$(0,-1)$};
\node at (-1,1) {$(-1,0)$};
\node at (0,-2.7) {$\underline{D(H)}$};

\end{tikzpicture}
& \hspace{0.2cm} &
\begin{tikzpicture}[scale=1.1]

\draw [->,line width=0.75mm] (0,0) coordinate (a) -- (0,2) coordinate (b1);
\draw [->,line width=0.75mm] (0,0) coordinate (a) -- (2,0) coordinate (b2);
\draw [->,line width=0.75mm] (0,0) coordinate (a1) -- (-1.7,-1.7) coordinate (b3);
\draw [->,line width=0.75mm] (0,0)-- (0,-2);
\draw [->,line width=0.75mm] (0,0)-- (-2,0);

\node at (2.7,0) {$ 0 $};
\node at (-2.3,0) {$1 $};
\node at (0,-2.3) {$1 $};
\node at (-2,-2) {$2 $};
\node at (0,2.2) {$0 $};

\node at (1,-1) {$(0,-1)$};
\node at (1,1) {$(0,0)$};
\node at (-1.7,-0.7) {$(-1,-1)$};

\node at (-0.8,-1.7) {$(-1,-1)$};
\node at (-1,1) {$(-1,0)$};
\node at (0,-2.7) {$\underline{D(2H-E_1-E_2)}$};

\end{tikzpicture} \end{array}$$
On the other hand, the monodromies $\gamma_{i}^j:=\gamma(z_{i} ; \sigma_j) \in G(z_{i},z_{0})$ are given as follows: $\gamma_{0}^k(\theta)=W$ for $k=1,...,5$ while

$$ \begin{array}{ccc} \gamma_{1}^1(\theta)= e^{-t} \cdot z_1+e^{-t} \cdot z_2+ \frac{1}{z_1}+\frac{e^{2 \pi i \theta}}{z_1 z_2} + \frac{e^{2 \pi i \theta} }{z_2} & ; & \gamma_{1}^2(\theta)= e^{-t} \cdot z_1+e^{-t} \cdot z_2+ \frac{1}{z_1}+\frac{e^{2 \pi i \theta}}{z_1 z_2} + \frac{e^{2 \pi i \theta}}{z_2} \end{array}$$

$$ \begin{array}{ccc} \gamma_{1}^3(\theta)= e^{-t} \cdot z_1+e^{-t+2 \pi i \theta} \cdot z_2+ \frac{1}{z_1}+\frac{1}{z_1 z_2} + \frac{1}{z_2} & ; & \gamma_{1}^4(\theta)= e^{-t} \cdot z_1+e^{-t+2 \pi i \theta} \cdot z_2+ \frac{1}{z_1}+\frac{1}{z_1 z_2} + \frac{1}{z_2} \end{array}$$

$$ \gamma_{1}^5(\theta)= e^{-t} \cdot z_1+e^{-t+2 \pi i \theta} \cdot z_2+ \frac{1}{z_1}+\frac{1}{z_1 z_2} + \frac{1}{z_2}$$

$$ \begin{array}{ccc} \gamma_{2}^1(\theta)= e^{-t} \cdot z_1+e^{-t} \cdot z_2+ \frac{e^{2 \pi i \theta}}{z_1}+\frac{e^{2 \pi i \theta}}{z_1 z_2} + \frac{1}{z_2} & ; & \gamma_{2}^2(\theta)= e^{-t+ 2 \pi i \theta} \cdot z_1+e^{-t} \cdot z_2+ \frac{1}{z_1}+\frac{1}{z_1 z_2} + \frac{1}{z_2} \end{array}$$

$$ \begin{array}{ccc} \gamma_{2}^3(\theta)= e^{-t+2 \pi i \theta} \cdot z_1+e^{-t} \cdot z_2+ \frac{1}{z_1}+\frac{1}{z_1 z_2} + \frac{1}{z_2} & ; & \gamma_{2}^4(\theta)= e^{-t+2 \pi i \theta} \cdot z_1+e^{-t} \cdot z_2+ \frac{1}{z_1}+\frac{1}{z_1 z_2} + \frac{1}{z_2} \end{array}$$

$$ \gamma_{2}^5(\theta)= e^{-t} \cdot z_1+e^{-t} \cdot z_2+ \frac{e^{2 \pi i \theta}}{z_1}+\frac{e^{2 \pi i \theta}}{z_1 z_2} + \frac{1}{z_2}$$

$$ \begin{array}{ccc} \gamma_{3}^1(\theta)= e^{-t} \cdot z_1+e^{-t} \cdot z_2+ \frac{e^{2 \pi i \theta}}{z_1}+\frac{e^{2 \pi i \theta}}{z_1 z_2} + \frac{e^{2 \pi i \theta}}{z_2} & ; & \gamma_{3}^2(\theta)= e^{-t+2 \pi i \theta} \cdot z_1+e^{-t} \cdot z_2+ \frac{1}{z_1}+\frac{1}{z_1 z_2} + \frac{e^{2 \pi i \theta}}{z_2} \end{array}$$

$$ \begin{array}{ccc} \gamma_{3}^3(\theta)= e^{-t+2 \pi i \theta} \cdot z_1+e^{-t} \cdot z_2+ \frac{1}{z_1}+\frac{1}{z_1 z_2} + \frac{e^{2 \pi i \theta}}{z_2} & ; & \gamma_{3}^4(\theta)= e^{-t} \cdot z_1+e^{-t+2 \pi i \theta} \cdot z_2+ \frac{e^{2 \pi i \theta}}{z_1}+\frac{1}{z_1 z_2} + \frac{1}{z_2} \end{array}$$

$$ \gamma_{3}^5(\theta)= e^{-t} \cdot z_1+e^{-t+2 \pi i \theta} \cdot z_2+ \frac{1}{z_1}+\frac{1}{z_1 z_2} + \frac{e^{2 \pi i \theta}}{z_2}$$

$$ \begin{array}{ccc} \gamma_{4}^1(\theta)= e^{-t} \cdot z_1+e^{-t} \cdot z_2+ \frac{e^{2 \pi i \theta}}{z_1}+\frac{e^{4 \pi i \theta}}{z_1 z_2} + \frac{e^{2 \pi i \theta}}{z_2} & ; & \gamma_{4}^2(\theta)= e^{-t+2 \pi i \theta} \cdot z_1+e^{-t} \cdot z_2+ \frac{e^{2 \pi i \theta}}{z_1}+\frac{1}{z_1 z_2} + \frac{e^{2 \pi i \theta}}{z_2} \end{array}$$

$$ \begin{array}{ccc} \gamma_{4}^3(\theta)= e^{-t+2 \pi i \theta} \cdot z_1+e^{-t+2 \pi i \theta} \cdot z_2+ \frac{1}{z_1}+\frac{1}{z_1 z_2} + \frac{1}{z_2} & ; & \gamma_{4}^4(\theta)= e^{-t+2 \pi i \theta} \cdot z_1+e^{-t+2 \pi i \theta} \cdot z_2+ \frac{1}{z_1}+\frac{1}{z_1 z_2} + \frac{1}{z_2} \end{array}$$

$$ \gamma_{4}^5(\theta)= e^{-t} \cdot z_1+e^{-t+2 \pi i \theta} \cdot z_2+ \frac{e^{2 \pi i \theta}}{z_1}+\frac{e^{2 \pi i \theta}}{z_1 z_2} + \frac{1}{z_2}.$$ The images of these paths under the Argument map, lifted to the universal cover of $\mathbb{T}^2 = \mathbb{R}^2/\mathbb{Z}^2$ are:

$$
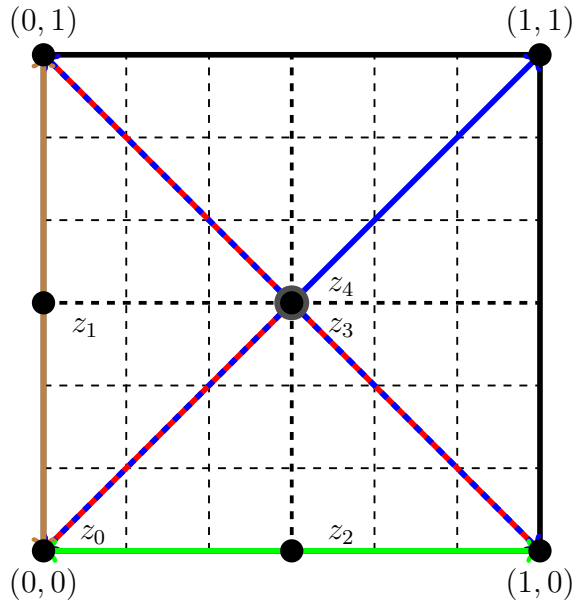
\begin{tikzpicture}[scale=1.1]

\draw [line width=0.75mm] (0,0)--(0,6);
\draw [line width=0.75mm] (0,0)--(6,0);
\draw [line width=0.75mm] (6,0)-- (6,6);
\draw [line width=0.75mm] (0,6)-- (6,6);

\draw [line width=0.25mm,dashed] (1,0)--(1,6);
\draw [line width=0.25mm,dashed] (2,0)-- (2,6);
\draw [line width=0.
5mm,dashed] (3,0)--(3,6);
\draw [line width=0.25mm,dashed] (4,0)-- (4,6);
\draw [line width=0.25mm,dashed] (5,0)--(5,6);

\draw [line width=0.25mm,dashed] (0,1)-- (6,1);
\draw [line width=0.25mm,dashed] (0,2)--(6,2);
\draw [line width=0.5mm,dashed] (0,3)-- (6,3);
\draw [line width=0.25mm,dashed] (0,4)--(6,4);
\draw [line width=0.25mm,dashed] (0,5)-- (6,5);

\draw [->,bend right,line width=0.75mm,red] (3,3)--(0,0);
\draw [->,bend right,line width=0.75mm,red] (3,3)--(0,6);
\draw [->,bend right, line width=0.75mm,red] (3,3)--(6,0);

\draw [->,line width=0.75mm,dashed,blue,bend left] (3,3)--(0,0);
\draw [->,bend left,line width=0.75mm,dashed,blue] (3,3)--(6,0);
\draw [->,out=150,in=30,line width=0.75mm,dashed,blue] (3,3)--(0,6);
\draw [->,line width=0.75mm,blue] (3,3)--(6,6);

\draw [->,line width=0.75mm,green] (3,0)--(6,0);
\draw [->,line width=0.75mm,green] (3,0)--(0,0);

\draw [->,line width=0.75mm,brown] (0,3)--(0,6);
\draw [->,line width=0.75mm,brown] (0,3)--(0,0);

\fill[black] (0,0) circle (4pt);
\fill[black] (0,3) circle (4pt);
\fill[black] (3,0) circle (4pt);
\fill[black!70] (3,3) circle (6pt);
\fill[black] (3,3) circle (4pt);
\fill[black] (0,6) circle (4pt);
\fill[black] (6,0) circle (4pt);
\fill[black] (6,6) circle (4pt);

\node at (3.6,0.2) {$z_{2}$};
\node at (0.5,2.7) {$z_{1}$};
\node at (3.6,2.7) {$z_{3}$};
\node at (3.6,3.2) {$z_4$};
\node at (0.6,0.2) {$z_{0}$};
\node at (0,-0.4) {$(0,0)$};
\node at (6,-0.4) {$(1,0)$};
\node at (0,6.4) {$(0,1)$};
\node at (6,6.4) {$(1,1)$};

\end{tikzpicture}$$

\hspace{-0.6cm} \underline{(4) Blow up of $\mathbb{P}^2$ at three points:} The fan $\Sigma$ for $X=Bl_3(\mathbb{P}^2)$ is given by
$$
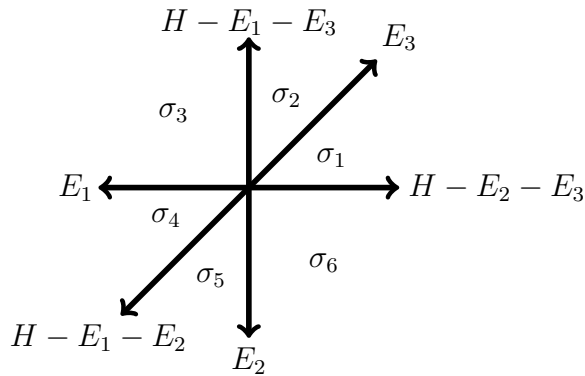
\begin{tikzpicture}

\draw [->,line width=0.75mm] (0,0) coordinate (a) -- (0,2) coordinate (b1);
\draw [->,line width=0.75mm] (0,0) coordinate (a) -- (2,0) coordinate (b2);
\draw [->,line width=0.75mm] (0,0) coordinate (a1) -- (-1.7,-1.7) coordinate (b3);
\draw [->,line width=0.75mm] (0,0)-- (0,-2);
\draw [->,line width=0.75mm] (0,0)-- (-2,0);
\draw [->,line width=0.75mm] (0,0)-- (1.7,1.7);

\node at (3.3,0) {$ H-E_2-E_3$};
\node at (0,2.2) {$H-E_1-E_3 $};
\node at (-2,-2) {$H-E_1-E_2 $};

\node at (-2.3,0) {$ E_1 $};
\node at (0,-2.3) {$E_2 $};
\node at (2,2) {$E_3$};

\node at (1.1,0.4) {$\sigma_1$};
\node at (0.5,1.2) {$\sigma_2$};
\node at (-1,1) {$\sigma_3$};
\node at (-1.1,-0.4) {$\sigma_4$};
\node at (-0.5,-1.2) {$\sigma_5$};
\node at (1,-1) {$\sigma_6$};

\end{tikzpicture}$$ In this case $\mathcal{E}=\left \{ 0 , H,H-E_1,H-E_2,H-E_3,2H-E_1-E_2-E_3 \right \} \subset Pic(Bl_3(\mathbb{P}^2))$. The divisors $D(E) \in Div_T(X)$ together with their corresponding weights are given as follows:

$$ \begin{array}{ccc}

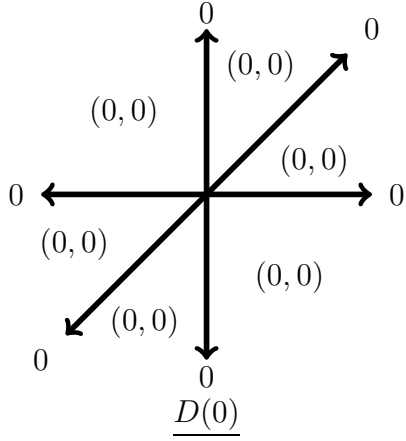
\begin{tikzpicture}[scale=1.1]

\draw [->,line width=0.75mm] (0,0) coordinate (a) -- (0,2) coordinate (b1);
\draw [->,line width=0.75mm] (0,0) coordinate (a) -- (2,0) coordinate (b2);
\draw [->,line width=0.75mm] (0,0) coordinate (a1) -- (-1.7,-1.7) coordinate (b3);
\draw [->,line width=0.75mm] (0,0)-- (0,-2);
\draw [->,line width=0.75mm] (0,0)-- (-2,0);
\draw [->,line width=0.75mm] (0,0)-- (1.7,1.7);

\node at (2.3,0) {$ 0 $};
\node at (0,2.2) {$0 $};
\node at (-2,-2) {$0 $};

\node at (-2.3,0) {$0 $};
\node at (0,-2.2) {$0 $};
\node at (2,2) {$0 $};

\node at (1.3,0.4) {$(0,0)$};
\node at (0.65,1.55) {$(0,0)$};
\node at (-1,1) {$(0,0)$};
\node at (-1.6,-0.6) {$(0,0)$};
\node at (-0.75,-1.55) {$(0,0)$};
\node at (1,-1) {$(0,0)$};
\node at (0,-2.7) {$\underline{D(0)}$};

\end{tikzpicture}
& \hspace{0.2cm} &
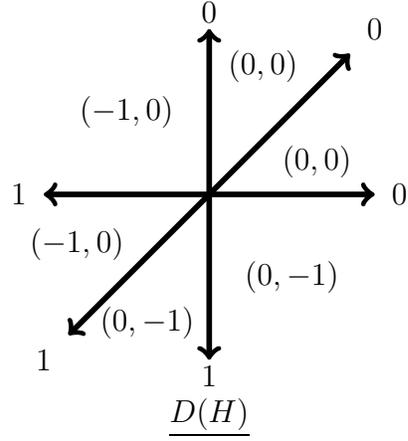
\begin{tikzpicture}[scale=1.1]

\draw [->,line width=0.75mm] (0,0) coordinate (a) -- (0,2) coordinate (b1);
\draw [->,line width=0.75mm] (0,0) coordinate (a) -- (2,0) coordinate (b2);
\draw [->,line width=0.75mm] (0,0) coordinate (a1) -- (-1.7,-1.7) coordinate (b3);
\draw [->,line width=0.75mm] (0,0)-- (0,-2);
\draw [->,line width=0.75mm] (0,0)-- (-2,0);
\draw [->,line width=0.75mm] (0,0)-- (1.7,1.7);

\node at (2.3,0) {$ 0 $};
\node at (0,2.2) {$0 $};
\node at (-2,-2) {$1 $};

\node at (-2.3,0) {$1 $};
\node at (0,-2.2) {$1 $};
\node at (2,2) {$0 $};

\node at (1.3,0.4) {$(0,0)$};
\node at (0.65,1.55) {$(0,0)$};
\node at (-1,1) {$(-1,0)$};
\node at (-1.6,-0.6) {$(-1,0)$};
\node at (-0.75,-1.55) {$(0,-1)$};
\node at (1,-1) {$(0,-1)$};

\node at (0,-2.7) {$\underline{D(H)}$};
\end{tikzpicture}
\end{array}$$
$$ \begin{array}{ccc}

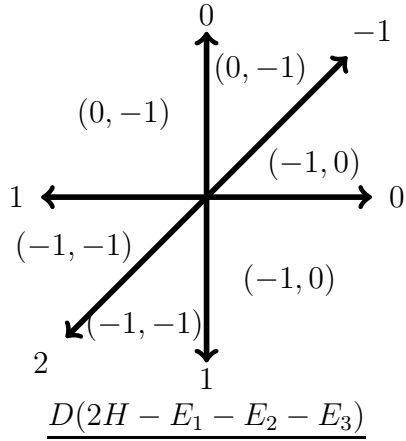
\begin{tikzpicture}[scale=1.1]

\draw [->,line width=0.75mm] (0,0) coordinate (a) -- (0,2) coordinate (b1);
\draw [->,line width=0.75mm] (0,0) coordinate (a) -- (2,0) coordinate (b2);
\draw [->,line width=0.75mm] (0,0) coordinate (a1) -- (-1.7,-1.7) coordinate (b3);
\draw [->,line width=0.75mm] (0,0)-- (0,-2);
\draw [->,line width=0.75mm] (0,0)-- (-2,0);
\draw [->,line width=0.75mm] (0,0)-- (1.7,1.7);

\node at (2.3,0) {$0 $};
\node at (0,2.2) {$0 $};
\node at (-2,-2) {$2 $};

\node at (-2.3,0) {$1 $};
\node at (0,-2.2) {$1 $};
\node at (2,2) {$-1 $};

\node at (1.3,0.4) {$(-1,0)$};
\node at (0.65,1.55) {$(0,-1)$};
\node at (-1,1) {$(0,-1)$};
\node at (-1.6,-0.6) {$(-1,-1)$};
\node at (-0.75,-1.55) {$(-1,-1)$};
\node at (1,-1) {$(-1,0)$};
\node at (0,-2.7) {$\underline{D(2H-E_1-E_2-E_3)}$};

\end{tikzpicture}
& \hspace{0.2cm} &
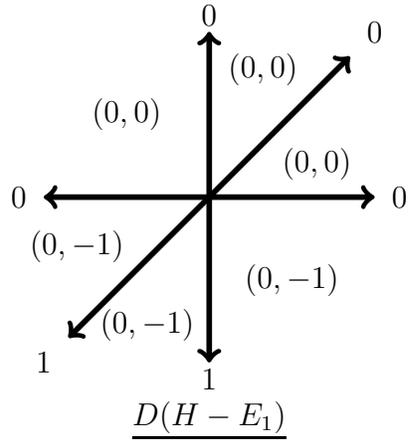
\begin{tikzpicture}[scale=1.1]

\draw [->,line width=0.75mm] (0,0) coordinate (a) -- (0,2) coordinate (b1);
\draw [->,line width=0.75mm] (0,0) coordinate (a) -- (2,0) coordinate (b2);
\draw [->,line width=0.75mm] (0,0) coordinate (a1) -- (-1.7,-1.7) coordinate (b3);
\draw [->,line width=0.75mm] (0,0)-- (0,-2);
\draw [->,line width=0.75mm] (0,0)-- (-2,0);
\draw [->,line width=0.75mm] (0,0)-- (1.7,1.7);

\node at (2.3,0) {$ 0 $};
\node at (0,2.2) {$0 $};
\node at (-2,-2) {$1$};

\node at (-2.3,0) {$0 $};
\node at (0,-2.2) {$1 $};
\node at (2,2) {$0 $};

\node at (1.3,0.4) {$(0,0)$};
\node at (0.65,1.55) {$(0,0)$};
\node at (-1,1) {$(0,0)$};
\node at (-1.6,-0.6) {$(0,-1)$};
\node at (-0.75,-1.55) {$(0,-1)$};
\node at (1,-1) {$(0,-1)$};

\node at (0,-2.7) {$\underline{D(H-E_1)}$};
\end{tikzpicture}
\end{array}$$

$$ \begin{array}{ccc}

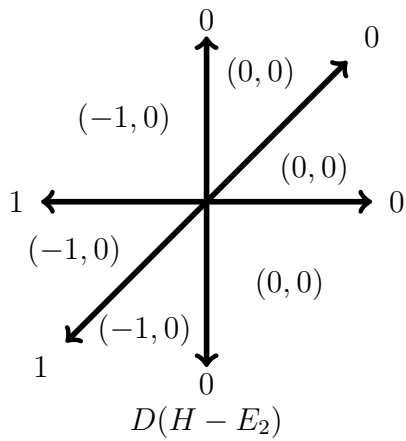
\begin{tikzpicture}[scale=1.1]

\draw [->,line width=0.75mm] (0,0) coordinate (a) -- (0,2) coordinate (b1);
\draw [->,line width=0.75mm] (0,0) coordinate (a) -- (2,0) coordinate (b2);
\draw [->,line width=0.75mm] (0,0) coordinate (a1) -- (-1.7,-1.7) coordinate (b3);
\draw [->,line width=0.75mm] (0,0)-- (0,-2);
\draw [->,line width=0.75mm] (0,0)-- (-2,0);
\draw [->,line width=0.75mm] (0,0)-- (1.7,1.7);

\node at (2.3,0) {$0 $};
\node at (0,2.2) {$0 $};
\node at (-2,-2) {$1 $};

\node at (-2.3,0) {$1 $};
\node at (0,-2.2) {$0 $};
\node at (2,2) {$0 $};

\node at (1.3,0.4) {$(0,0)$};
\node at (0.65,1.55) {$(0,0)$};
\node at (-1,1) {$(-1,0)$};
\node at (-1.6,-0.6) {$(-1,0)$};
\node at (-0.75,-1.55) {$(-1,0)$};
\node at (1,-1) {$(0,0)$};

\node at (0,-2.7) {$\underline{D(H-E_2)}$};

\end{tikzpicture}
& \hspace{0.2cm} &
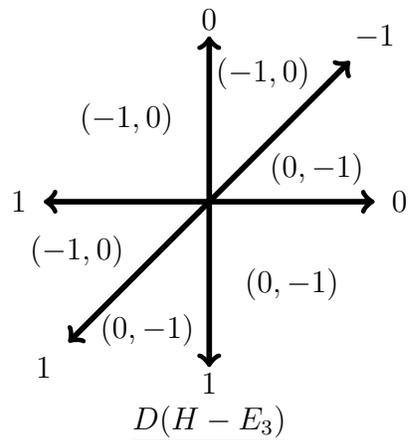
\begin{tikzpicture}[scale=1.1]

\draw [->,line width=0.75mm] (0,0) coordinate (a) -- (0,2) coordinate (b1);
\draw [->,line width=0.75mm] (0,0) coordinate (a) -- (2,0) coordinate (b2);
\draw [->,line width=0.75mm] (0,0) coordinate (a1) -- (-1.7,-1.7) coordinate (b3);
\draw [->,line width=0.75mm] (0,0)-- (0,-2);
\draw [->,line width=0.75mm] (0,0)-- (-2,0);
\draw [->,line width=0.75mm] (0,0)-- (1.7,1.7);

\node at (2.3,0) {$ 0 $};
\node at (0,2.2) {$0 $};
\node at (-2,-2) {$1 $};

\node at (-2.3,0) {$1 $};
\node at (0,-2.2) {$1 $};
\node at (2,2) {$-1 $};

\node at (1.3,0.4) {$(0,-1)$};
\node at (0.65,1.55) {$(-1,0)$};
\node at (-1,1) {$(-1,0)$};
\node at (-1.6,-0.6) {$(-1,0)$};
\node at (-0.75,-1.55) {$(0,-1)$};
\node at (1,-1) {$(0,-1)$};

\node at (0,-2.7) {$\underline{D(H-E_3)}$};
\end{tikzpicture}
\end{array}$$ On the other hand, the monodromies $\gamma_{i}^j:=\gamma(z_{i} ; \sigma_j) \in G(z_{i},z_{0})$ are given as follows: $\gamma_{0}^k(\theta)=W$ for $k=1,...,6$ while

$$ \gamma_{1}^1(\theta)= \gamma_1^2 (\theta)= z_1+ z_2+z_1 z_2+ \frac{e^{2 \pi i \theta}}{z_1}+\frac{e^{2 \pi i \theta}}{z_1 z_2} + \frac{e^{2 \pi i \theta} }{z_2} $$

$$ \gamma_{1}^3(\theta)= \gamma_1^4( \theta)= e^{2 \pi i \theta} \cdot z_1+ z_2+e^{ 2 \pi i \theta} \cdot z_1 z_2+ \frac{1}{z_1}+\frac{1}{z_1 z_2} + \frac{e^{2 \pi i \theta} }{z_2} $$

$$ \gamma_{1}^5(\theta)= \gamma_1^6( \theta)= z_1+ e^{2 \pi i \theta} \cdot z_2+e^{ 2 \pi i \theta} \cdot z_1 z_2+ \frac{e^{2 \pi i \theta}}{z_1}+\frac{1}{z_1 z_2} + \frac{1 }{z_2} $$

$$ \gamma_{2}^1(\theta)= \gamma_2^6 (\theta)= z_1+ e^{ 2 \pi i \theta} \cdot z_2+z_1 z_2+ \frac{e^{2 \pi i \theta}}{z_1}+\frac{e^{2 \pi i \theta}}{z_1 z_2} + \frac{1 }{z_2} $$

$$ \gamma_{2}^2(\theta)= \gamma_2^3 (\theta)= e^{ 2 \pi i \theta} \cdot z_1+ z_2+z_1 z_2+ \frac{1}{z_1}+\frac{e^{2 \pi i \theta}}{z_1 z_2} + \frac{e^{2 \pi i \theta} }{z_2} $$

$$ \gamma_{2}^4(\theta)= \gamma_2^5 (\theta)= e^{ 2 \pi i \theta} \cdot z_1+ e^{2 \pi i \theta} \cdot z_2+e^{2 \pi i \theta} \cdot z_1 z_2+ \frac{1}{z_1}+\frac{1}{z_1 z_2} + \frac{1 }{z_2} $$

$$ \gamma_{3}^1(\theta)= \gamma_3^2 (\theta)= \gamma_3^3(\theta) = z_1+ z_2+z_1 z_2+ \frac{1}{z_1}+\frac{e^{2 \pi i \theta}}{z_1 z_2} + \frac{e^{2 \pi i \theta} }{z_2} $$

$$ \gamma_{3}^4(\theta)= \gamma_3^5 (\theta)= \gamma_3^6(\theta) = z_1+ e^{2 \pi i \theta} \cdot z_2+e^{2 \pi i \theta} \cdot z_1 z_2+ \frac{1}{z_1}+\frac{1}{z_1 z_2} + \frac{1 }{z_2} $$

$$ \gamma_{4}^1(\theta)= \gamma_4^2 (\theta)= \gamma_4^6(\theta) = z_1+ z_2+z_1 z_2+ \frac{e^{2 \pi i \theta}}{z_1}+\frac{e^{2 \pi i \theta}}{z_1 z_2} + \frac{1 }{z_2} $$

$$ \gamma_{4}^3(\theta)= \gamma_4^4 (\theta)= \gamma_4^5(\theta) = e^{2 \pi i \theta} \cdot z_1+ z_2+e^{2 \pi i \theta} \cdot z_1 z_2+ \frac{1}{z_1}+\frac{1}{z_1 z_2} + \frac{1 }{z_2} $$

$$ \gamma_{5}^1(\theta)= \gamma_5^5 (\theta)= \gamma_5^6(\theta) = z_1+ e^{ 2 \pi i \theta} \cdot z_2+z_1 z_2+ \frac{e^{2 \pi i \theta}}{z_1}+\frac{1}{z_1 z_2} + \frac{1 }{z_2} $$

$$ \gamma_{5}^2(\theta)= \gamma_5^3 (\theta)= \gamma_5^4(\theta) = e^{2 \pi i \theta} \cdot z_1+ z_2+ z_1 z_2+ \frac{1}{z_1}+\frac{1}{z_1 z_2} + \frac{e^{2 \pi i \theta} }{z_2}. $$ The images of these paths under the Argument map, lifted to the universal cover of $\mathbb{T}^2 = \mathbb{R}^2/\mathbb{Z}^2$ are:

$$
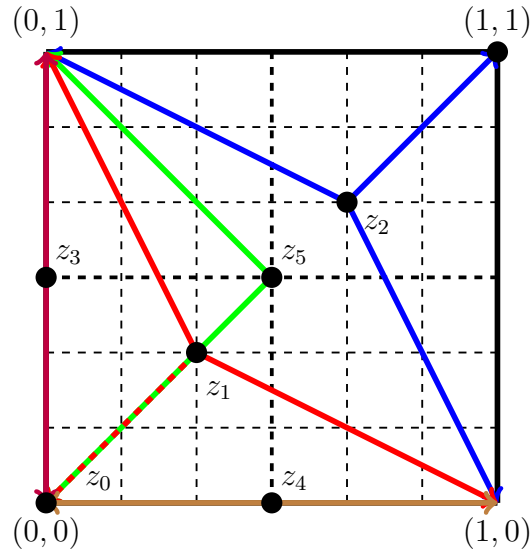
\begin{tikzpicture}

\draw [line width=0.75mm] (0,0)--(0,6);
\draw [line width=0.75mm] (0,0)--(6,0);
\draw [line width=0.75mm] (6,0)-- (6,6);
\draw [line width=0.75mm] (0,6)-- (6,6);

\draw [line width=0.25mm,dashed] (0,1)--(6,1);
\draw [line width=0.25mm,dashed] (0,2)--(6,2);
\draw [line width=0.5mm,dashed] (0,3)--(6,3);
\draw [line width=0.25mm,dashed] (0,4)--(6,4);
\draw [line width=0.25mm,dashed] (0,5)--(6,5);

\draw [line width=0.25mm,dashed] (1,0)--(1,6);
\draw [line width=0.25mm,dashed] (2,0)--(2,6);
\draw [line width=0.5mm,dashed] (3,0)--(3,6);
\draw [line width=0.25mm,dashed] (4,0)--(4,6);
\draw [line width=0.25mm,dashed] (5,0)--(5,6);

\draw [->,line width=0.75mm,blue] (4,4)--(6,6);
\draw [->,line width=0.75mm,blue] (4,4)--(6,0);
\draw [->,line width=0.75mm,blue] (4,4)--(0,6);

\draw [->,line width=0.75mm,green] (3,3)--(0,0);
\draw [->,line width=0.75mm,green] (3,3)--(0,6);

\draw [->,line width=0.75mm,dashed,red] (2,2)--(0,0);
\draw [->,line width=0.75mm,red] (2,2)--(6,0);
\draw [->,line width=0.75mm,red] (2,2)--(0,6);

\draw [->,line width=0.75mm,purple] (0,3)--(0,6);
\draw [->,line width=0.75mm,purple] (0,3)--(0,0);

\draw [->,line width=0.75mm,brown] (3,0)--(6,0);
\draw [->,line width=0.75mm,brown] (3,0)--(0,0);

\fill[black] (0,0) circle (4pt);
\fill[black] (2,2) circle (4pt);
\fill[black] (4,4) circle (4pt);
\fill[black] (0,3) circle (4pt);
\fill[black] (3,0) circle (4pt);
\fill[black] (3,3) circle (4pt);
\fill[black] (6,6) circle (4pt);

\node at (0.7,0.3) {$z_0$};
\node at (3.3,0.3) {$z_4$};
\node at (0.3,3.3) {$z_3$};
\node at (3.3,3.3) {$z_5$};
\node at (4.4,3.7) {$z_2$};
\node at (2.3,1.5) {$z_1$};

\node at (0,-0.4) {$(0,0)$};
\node at (6,-0.4) {$(1,0)$};
\node at (0,6.4) {$(0,1)$};
\node at (6,6.4) {$(1,1)$};
\end{tikzpicture}$$
$\square$

\bigskip

\hspace{-0.6cm} \bf Acknowledgements: \rm The work was partially supported by ERC grant AdG 669655.

\end{document}